\UseRawInputEncoding
\RequirePackage{fix-cm}
\documentclass[smallextended]{svjour3}       
\smartqed  
\usepackage{graphicx}
\usepackage{color,ragged2e,newlfont,pifont,stmaryrd,xspace}
\usepackage{amsmath,amssymb,amsfonts,textcomp,yhmath}
\usepackage{cases}
\usepackage{dsfont}
\usepackage{caption}
\usepackage{subfigure}
\newenvironment{prooof}[1][\underline{Proof}]
{\noindent\textbf{#1:} }{\ \rule{0.5em}{0.5em}}
\begin{document}

\title{Global analysis of a spatiotemporal cellular model for the transmission of hepatitis C virus with Hattaf-Yousfi functional response
}

\titlerunning{Spatiotemporal cellular model for the HCV infection} 

\author{Alexis Nangue \and Bruno Nde Tchiffo
}


\institute{University of Maroua, Higher Teachers' Training College,
Department of Mathematics \at
               P.O.Box 55 Maroua, Cameroon \\
              Tel.: +237 696 419 142\\
              \email{alexnanga02@yahoo.fr}           
           \and
           University of Maroua, Faculty of
Science, Department of Mathematics and Computer Science, Cameroon
\at
              P.O. Box 814 Maroua\\
              \email{ndebruno2@gmail.com}
}

\date{Received: date / Accepted: date}

\maketitle

\begin{abstract}
In this paper, a mathematical analysis of the global dynamics of a
partial differential equation viral infection cellular model  is
carried out. We study the dynamics of a hepatitis C virus (HCV)
model, under therapy, that considers both absorption phenomenon and diffusion of virions, infected and uninfected cells in
liver. Firstly, we prove boundedness of
the potential solutions, global existence, uniqueness and positivity
of the solutions  to the obtained initial value and boundary
problem. Then, the dynamical behavior of the model is completely
determined by a threshold parameter called the basic reproduction
number $\mathcal{R}_{0}$. We show that the uninfected spatially
homogeneous equilibrium of the model is globally asymptotically
stable if $\mathcal{R}_{0} \leq 1 $ by using the direct Lyapunov
method. This means that the HCV is cleared and the disease dies out.
Also, the global asymptotically properties stability of the infected
spatially homogeneous equilibrium of the model are studied via a
skilful construction of a suitable Lyapunov functional. It means
that the HCV persists in the host and the infection becomes chronic.
Finally, numerical simulations are performed to support the
theoretical results obtained.

\end{abstract}
\section{Introduction}\label{intro}
The dynamics of viruses, in particular the dynamics of the hepatitis
C virus, remains a very active field of research in the world of
science. Moreover, the 2020 Nobel Prize in Medicine was awarded to
three researchers, namely the British Michael Hougton and the
Americans Harvey Alter and Charles Rice. They were awarded this
Nobel Prize for their very advanced research work on the hepatitis C
virus.

According to World Health Organization(WHO) \cite{who1}, 71 million persons were living with
chronic hepatitis C virus (HCV) infection worldwide and  399 000
persons had died from cirrhosis or hepatocellular carcinoma
following a survey done in 2015. Aside from the burden of HCV
infection secondary to liver-related sequelae, HCV causes an
additional burden through comorbidities among persons with HCV
infection, including depression, diabetes mellitus and chronic renal
disease. In May 2016, the World Health Assembly endorsed the Global
Health Sector Strategy for 2016-2021 on viral hepatitis (HBV and HCV
infection), which proposes to eliminate viral hepatitis as a public
health threat by 2030. Elimination is defined as a 90\% reduction in
new chronic infections and a 65\% reduction in mortality compared
with the 2015 baseline.

 Mathematicians cannot stay aside from
disastrous situation decried by WHO. In view of the vital importance
of the liver and the aforementioned facts, any contribution to a
better understanding of HCV infection process is of great interest.
Mathematical models have been developed to help understand and
control the dynamics of HCV within an infected host such as in
\cite{chong2015,Rongetal2013,chatterjeeetal2012,Guedjetneumann2010}.
 The dynamics of viral infections such
as the Ebola virus disease(EVD), the human immunodeficiency virus
(HIV) infection, the hepatitis B virus (HBV) infection, the
hepatitis C virus (HCV) infection and, new corona virus infection
have been modeled mathematically in a host. One of the earliest
temporal models was the within-host basic viral infection model
proposed in \cite{nowakbagham} to study HIV infection, and later
adopted to HBV \cite{nowakbonofer,ciupeetal2007}. Particularly,
numerous mathematical models describing the temporal dynamics of HCV
have been initially proposed by Neumann and al \cite{neumanetal1998}
using the classical viral infection cellular model, and later have
been extended in
\cite{Daharietal2007,Guedjetneumann2010,chatterjeeetal2012,Rongetal2013}.
Motivated by what has been done in
\cite{nowakbonofer,ciupeetal2007,neumanetal1998}, Chong and
al.\cite{chong2015} formulated the basic HCV temporal intra-host
model or cellular model infection with therapy as a system of three
differential equations :
\begin{numcases}\strut
\nonumber \frac{d H(t)}{d t}=\lambda-d H(t)-(1-\eta)\beta H(t)V(t),\\
          \frac{d I(t)}{d t}=(1-\eta)\beta H(t) V(t)-\alpha I(t), \label{basicmodthe}\\
\nonumber \frac{d V(t)}{d t}=(1-\varepsilon)k I(t)-\mu V(t),
\end{numcases}
where the equations relate the dynamics relationship between, H as
the uninfected target cells (hepatocytes), I as the infected cells
and V as the viral load (amount of viruses present in the blood).
 In the system (\ref{basicmodthe})
the key assumption is that hepatocytes and viruses are well mixed,
and ignores the mobility of hepatocytes C viruses, the infected and uninfected target cells. To study the
influences of spatial structures of virus dynamics, Wang and Wang in
\cite{wangwang} assuming that the motion of virus follows Fickian
diffusion, that is to say, the population flux of virus is
proportional to the concentration gradient and the proportionality
constant is taken to be negative \cite{fichian}.
\\\noindent Moreover, in system (\ref{basicmodthe}), the rate of
infection is assumed to be bilinear in the virus V and uninfected
hepatocytes T. It is shown in \cite{minetal} that this bilinear rate
of infection could be unrealistic. However, the actual incidence
rate is probably not linear over the entire range of T and V. Thus
is reasonable to assume that the infection rate is given by a more
general one, known as the Hattaf-Yousfi functionnal response
\cite{HattafK2016} of the form $\frac{\beta H V}{\alpha_{0}
   +\alpha_{1}H+\alpha_{2}V+\alpha_{3}HV}$ where $\alpha_{0}>0$, $\alpha_{1}\geq0$,
    $\alpha_{2}\geq0$, $\alpha_{3}\geq0$ are constants.
The function
    $\frac{\beta H
}{\alpha_{0}
   +\alpha_{1}H+\alpha_{2}V+\alpha_{3}HV}$ satisfies the hypotheses
   $(H_{1})$, $(H_{2})$ and $(H_{3})$ of general incidence rate
   presented in
   \cite{Hattafetal2012,Hattafetal2013a,Hattafetal2013b,HY1}. The
   Hattaf-Yousfi type of functional response was introduced by
   Hattaf and al. in \cite{HattafK2016}. This functional response
   generalizes many functional responses and it was used in
   \cite{riad2016} to describe the dynamics of labour market. Thus,
    when $\alpha_{0}= 1$,
    the Hattaf-Yousfi functional response is reduced to
the specific functional response used by Hattaf and al in
\cite{Hattafetal2015}. Furthermore, if
$\alpha_{3}=\alpha_{1}\alpha_{2}$ and $\alpha_{0}= 1$, the
Hattaf-Yousfi functional response is reduced to Crowley-Martin
functional response \cite{crowleym1989} and was used in
\cite{zhouandcui2011}. When $\alpha_{3}=0$ et $\alpha_{0}=1$ the
Hattaf-Yousfi functional response is simplified to
Beddington-DeAngelis functional
response\cite{Bedding1975,Deangelisetal1975}, and was used in
\cite{huang2009,huang2011,wangetal2011,zhangxu2013}. When
$\alpha_{1}>0$, $\alpha_{2}=\alpha_{3}=0$ and $\alpha_{0}=1$, the
Hattaf-Yousfi functional response is reduced to Holling type II
functional response\cite{Liandma2007}. And when
$\alpha_{1}=\alpha_{3}=0$, $\alpha_{2}>0$ and $\alpha_{0}=1$ it
expresses a saturation response\cite{Songneuman2007}. Moreover, when
$\alpha_{1}=\alpha_{2}=\alpha_{3}=0$, and  $\alpha_{0}=1$ the
Hattaf-Yousfi functional response is reduced to the mass action
principle(or Holling type I functional response). Also ordinary
differential system (\ref{basicmodthe})take into consideration the
cure of infected hepatocytes.
\\\indent
In this work, motivated by the breaches observed in the analysis and
the formulation of system (\ref{basicmodthe}), we construct and
analyze a PDE-cellular model system for HCV infection, which derives
from system (\ref{basicmodthe}) by incorporating the space,
Hattaf-Yousfi incidence rate, absorption effect and spontaneous
cure. It is worth mentioning that in \cite{chong2015} the authors
used mass-action kinetics for viral infection, neglected the cure
rate, ignored the absorption effect and the diffusion of free virions, 
susceptible cells and infected cells.
So the obtained model is an extension of the one in the first of
work done by  Chong and al.\cite{chong2015}.
\\\indent
The work is organized as follows. In section~\ref{sec:1}, we model
the phenomenon described through a system of partial differential
equations which leads to a initial value and boundary problem.
Section 3 is devoted to the study of the existence and uniqueness of
the global solution of our initial and boundary value problem, and
of the properties of this solution, namely positivity and
boundedness. Section~\ref{sec:3} deals with the stability and the
analysis of spatially homogeneous equilibria and numerical
simulations in section~\ref{sec:4}. We conclude our work and provide
a discussion in section~\ref{sec:5}.

\section{Formulation of the PDE-cellular model}\label{sec:1}
Let  $\Omega \subset \mathds{R}^{3} $ be the domain representing the
liver. Let $t \geq 0 $ be a given time and
$x=(x_{1}, x_{2},
x_{3})
 \in \Omega$. Denote respectively by $H(x, t)$,
 $I(x, t)$ and $V(x, t)$ the
concentrations of healthy hepatocytes, HCV infected hepatocytes, and
free HCV virions at time $t$ and location $x$. The dynamics of HCV
infection intra-host is the result of the dynamics
of each
compartment H, I, and V, and the various interactions
 between them.
We now describe the evolution of each compartment.
\paragraph{Fluctuation of healthy hepatocytes.}
Let $\nu$ be an elementary volume in $\Omega$. The variation of the
quantity of healthy hepatocytes in $\nu$ is described
under the
following assumptions. Healthy hepatocytes are produced
 at constant rate $\lambda$
from the bone marrow and die at rate $dH$. Virions
infect the
healthy hepatocytes at the rate
 $\frac{\beta H V}{\alpha_{0}
   +\alpha_{1}H+\alpha_{2}V+\alpha_{3}HV},$ where
 $\beta$ is the rate of transmission of the infection
 and $\alpha_{j}$,
     $j=0,1,2,3$ are positive constants.
 This generalized
incidence function replaces the mass-action function
 which has been
shown to cause unrealistic conditions for
 successful chronic HCV
infection. $\rho I$ is the cure rate of infected hepatocytes either
by noncytolytic mechanism or immunity or treatment. In addition, the
therapeutic effect of
 treatment in this model
involved the reduction of new infections, which
 is described in a fraction as $(1-\eta)$. The spatial
 motion of healthy hepatocytes follows the
Fickian diffusion law. Thus,
 the variation of healthy hepatocytes
 is expressed by the following
equation:
\begin{eqnarray*}
  \frac{\partial H}{\partial t}&=&D_{1}\Delta H(x,t)+\lambda-dH(x,t)\\
 &&-\frac{(1-\eta)\beta H(x,t)V(x,t)}{\alpha_{0}
  +\alpha_{1}H(x,t)+\alpha_{2}V(x,t)+\alpha_{3}H(x,t)V(x,t)} +\rho I(x,t),
\end{eqnarray*}

where $D_{1}$ represents the Healthy hepatocytes diffusion coefficient and
$$\Delta=\frac{\partial^{2}}{\partial x_{1}^{2}}
  + \frac{\partial^{2}}{\partial x_{2}^{2}}
  + \frac{\partial^{2}}{\partial x_{3}^{2}}$$
is the usual Laplacian operator in three-dimensional
 space.
\paragraph{Fluctuation of HCV infected cells.}
 The HCV infected cells die at rate $\alpha$
 per day so that $\frac{1}{\alpha}$ is
the life-expectancy of HCV infected hepatocytes. Healthy hepatocytes
become infected at the rate $\frac{\beta H V}{\alpha_{0}
   +\alpha_{1}H+\alpha_{2}V+\alpha_{3}HV}$. The spatial
 motion of HCV infected cells follows the
Fickian diffusion law. Thus, the
variation of infected hepatocytes is expressed by the following equation
\begin{equation*}
   \frac{\partial I}{\partial t}=D_{2}\Delta I(x,t)+\frac{(1-\eta)\beta H(x,t)V(x,t)}
   {\alpha_{0}+\alpha_{1}H(x,t)+\alpha_{2}V(x,t)+\alpha_{3}H(x,t)V(x,t)}
   -(\alpha+\rho)I(x,t),
 \end{equation*}
 where $D_{2}$ represents the HCV infected cells diffusion coefficient.
\paragraph{Fluctuation of free HCV virions.}
The infected hepatocytes produce virus at rate $kI$, and virus is
cleared at the rate $\mu V$. Also, the population of virions
decreases due to the infection at the rate $\frac{u (1-\eta)\beta
HV}
   {\alpha_{0}+\alpha_{1}H+\alpha_{2}V+\alpha_{3}HV}$
   due to absorption effect, where $u \in¸\{0,1\}$. The spatial
 motion of virions follows the
Fickian diffusion law.  In addition,
the therapeutic effect of
 treatment in this model
involved blocking virions production
(referred to as drug
effectiveness) which, is described
in fraction
 $(1-\varepsilon)$.
Thus, the variation of free virions is expressed by the following equation:
\begin{eqnarray*}
   \frac{\partial V}{\partial t} &=&D_{3}\Delta V(x,t)+
   (1-\varepsilon)kI(x,t)-\mu V(x,t)  \\
  &&-u\frac{(1-\eta)\beta H(x,t)V(x,t)}
   {\alpha_{0}+\alpha_{1}H(x,t)+\alpha_{2}V(x,t)+\alpha_{3}H(x,t)V(x,t)},
\end{eqnarray*}
where $D_{3}$ represents the free HCV virions diffusion coefficient.
\paragraph{The initial boundary value problem associated to PDE-cellular
model.}
In this section, we use the previous equations describing variables
variations to set up a complete PDE system modeling biological
dynamics for HCV infection. Let $T > 0 $ be a fixed time and define
\begin{equation*}
    \Omega_{T}=\Omega\times (0,T).
\end{equation*}
Therefore, in $\Omega_{T}$ the full system of PDE governing the HCV
infection becomes :
\begin{numcases}\strut
\nonumber \frac{\partial H}{\partial t}=D_{1}\Delta H(x,t)+\lambda-dH-\frac{(1-\eta)\beta HV}
{\alpha_{0}+\alpha_{1}H+\alpha_{2}V+\alpha_{3}HV}+\rho I,\\
          \frac{\partial I}{\partial t}=D_{2}\Delta I(x,t)+\frac{(1-\eta)\beta HV}{\alpha_{0}
          +\alpha_{1}H+\alpha_{2}V+\alpha_{3}HV}-(\alpha+\rho)I,\label{equation1.0}\\
\nonumber \frac{\partial V}{\partial t}=D_{3}\Delta
V+(1-\varepsilon)kI-\mu V-\frac{u(1-\eta)\beta
HV}{\alpha_{0}+\alpha_{1}H+\alpha_{2}V+\alpha_{3}HV}.
\end{numcases}

We use the Neumann homogeneous boundary conditions:
\begin{equation}\label{neuman}
\frac{\partial H}{\partial\eta}=\frac{\partial I}{\partial\eta}=\frac{\partial V}{\partial\eta}=0 \;\; \mbox{on} \;
\partial\Omega\times [0,T],
\end{equation}
where   $\frac{\partial }{\partial\eta}$ denotes the outward normal
derivative on $\partial\Omega$. The initial conditions are the
following :
\begin{equation}\label{initcond}
 H(x,0)=H_{0},\;\; I(x,0)=I_{0},\;\; V(x,0)=V_{0},\;\; x\in \Omega.
\end{equation}
The boundary conditions in (\ref{neuman}) imply that the Healthy hepatocytes, the HCV infected cells and free HCV
virions do not move across the boundary $\partial\Omega$. For an
epidemiological significance, we assume that the initial conditions
are positive and H\"{o}lder continuous, and satisfy $\frac{\partial
H_{0}}{\partial\eta}=\frac{\partial
I_{0}}{\partial\eta}=\frac{\partial
V_{0}}{\partial\eta}=0$ on
$\partial\Omega$.\\
\indent We then obtain the following initial boundary value problem
(IBVP) associated to our PDE-cellular model:
\begin{numcases}\strut
\nonumber \frac{\partial H}{\partial
t}=D_{1}\Delta H+\lambda-dH-\frac{(1-\eta)\beta HV}
{\alpha_{0}+\alpha_{1}H+\alpha_{2}V+\alpha_{3}HV}+\rho I\;\; \mbox{in}\;\; \Omega_{T},\\
\nonumber \frac{\partial I}{\partial t}=D_{2}\Delta I+\frac{(1-\eta)\beta
HV}{\alpha_{0}
          +\alpha_{1}H+\alpha_{2}V+\alpha_{3}HV}-(\alpha+\rho)I
          \;\; \mbox{in}\;\; \Omega_{T},\\
 \frac{\partial V}{\partial t}=D_{3}\Delta
V+(1-\varepsilon)kI-\mu V-\frac{u(1-\eta)\beta
HV}{\alpha_{0}+\alpha_{1}H+\alpha_{2}V+\alpha_{3}HV}\;\;
\mbox{in}\;\; \Omega_{T},\;\; \label{equation1.1}\\
 \nonumber
\frac{\partial H}{\partial\eta}=\frac{\partial I}{\partial\eta}=\frac{\partial V}{\partial\eta}=0  \;\;
\mbox{on} \;\; \partial\Omega\times [0,T], \\
\nonumber    H(x,0)=H_{0},\; I(x,0)=I_{0},\; V(x,0)=V_{0},\;\;  x
\in \Omega,
\end{numcases}
on which our study will focus on.
\section{Quantitative analysis and some properties of the
solutions for IBVP (\ref{equation1.1}) }\label{sec:2} In this
section, we provide a thorough study of the dynamics of IBVP
(\ref{equation1.1}) which yields various outcomes. Precisely, we
prove existence, uniqueness, positivity and boundedness of solutions
for IBVP (\ref{equation1.1}). This is done by combining variational
method and semigroups techniques to some useful functional analysis
arguments.
\subsection{Local existence and uniqueness of solutions
for the IBVP (\ref{equation1.1})}
Set
\begin{equation}\label{equation2.1}
F(H,I,V)=(F_{1}(H,I,V), F_{2} (H,I,V), F_{3}(H,I,V))^{T}
\end{equation}
where
\begin{equation*}
    F_{1}(H,I,V)=\lambda-dH-\frac{(1-\eta)
\beta HV}{\alpha_{0}+\alpha_{1}H+\alpha_{2}V+\alpha_{3}HV}+\rho I,
\end{equation*}
\begin{equation*}
    F_{2}(H,I,V)=\frac{(1-\eta)\beta HV}{\alpha_{0}
+\alpha_{1}H+\alpha_{2}V+\alpha_{3}HV}-(\alpha+\rho)I,
\end{equation*}
and
\begin{equation*}
 F_{3}(H,I,V)=(1-\varepsilon)kI-\mu V-u\frac{(1
-\eta)\beta HV}{\alpha_{0}+\alpha_{1}H+\alpha_{2}V+\alpha_{3}HV}
.
\end{equation*}
We have the following result which guarantees that the right-hand
side, without diffusion, of the PDE-model system (\ref{equation1.1})
is Lipschitz.
\begin{proposition}\label{proposition2.1.1}
Let $T\in \mathds{R}_{+}^{*}$ and $(H,I,V)\in
(C_{B}^{0}(\Omega\times[0,T]))^{3}$,
 where $C_{B}^{0}(\Omega\times[0,T])$
is the space of bounded and continuous functions on
 $\Omega\times[0,T]$. We suppose that F in (\ref{equation2.1})
  is defined on $L^{2}(\Omega\times (0,T))$.
Then $F_{1}, F_{2}$ and $F_{3}$ are uniformly Lipschitz continuous
on $L^{2}(\Omega\times(0,T))$ with respect to H, I and V.
\end{proposition}
\begin{prooof} Let $T\in \mathbb{R}_{+}^{*}$ and $(H_{1},I_{1},V_{1}),
(H_{2},I_{2},V_{2})\in (C_{B}^{0}(\Omega\times[0,T]))^{3}$. First,
by direct computation, we have :
\begin{equation}\label{lipF1}
\|F_{1}(H_{1},I_{1},V_{1})-F_{1}(H_{2},I_{2},V_{2})\|_{2} \leq
K_{1}^{1}\|H_{1}-H_{2}\|_{2}+K_{2}^{1}\|I_{1}-I_{2}\|_{2}
+K_{3}^{1}\|V_{1}-V_{2}\|_{2},
\end{equation}
with
\begin{equation}\label{coefF1}
K_{1}^{1}=d+(1-\eta)\beta\left(\frac{1}{\alpha_{2}}
+\frac{V_{m}}{\alpha_{0}}\right), \;\; K_{2}^{1}=\rho, \;\;
\;\; K_{3}^{1}=(1-\eta)\beta
\left(\frac{1}{\alpha_{1}}+\frac{H_{m}}{\alpha_{0}}\right),
\end{equation}
with  $H_{m}$ and $V_{m}$ given below.\\
Then
\begin{equation}\label{lipF2}
\|F_{2}(H_{1},I_{1},V_{1})-F_{2}(H_{2},I_{2},V_{2})\|_{2} \leq
K_{1}^{2}\|H_{1}-H_{2}\|_{2}+K_{2}^{2}\|I_{1}-I_{2}\|_{2}
+K_{3}^{2}\|V_{1}-V_{2}\|_{2},
\end{equation}
with
\begin{equation}\label{coefF2}
K_{1}^{2}=(1-\eta)\beta\left(\frac{1}{\alpha_{2}}
+\frac{V_{m}}{\alpha_{0}}\right), \;\; K_{2}^{2}=(\alpha+\rho)\;\;
and \;\; K_{3}^{2}=(1-\eta)\beta\left(\frac{1}{\alpha_{1}}
+\frac{H_{m}}{\alpha_{0}}\right).
\end{equation}
Finally
\begin{equation}\label{lipF3}
\|F_{3}(H_{1},I_{1},V_{1})-F_{3}(H_{2},I_{2},V_{2})\|_{2} \leq
K_{1}^{3}\|H_{1}-H_{2}\|_{2}+K_{2}^{3}\|I_{1}-I_{2}\|_{2}
+K_{3}^{3}\|V_{1}-V_{2}\|_{2},
\end{equation}
with
\begin{equation}\label{coefF3}
K_{1}^{3}=u(1-\eta)\beta
\left(\frac{1}{\alpha_{2}}+\frac{V_{m}}{\alpha_{0}}\right), \;\;
K_{2}^{3}=k(1-\varepsilon), \;\; and \;\; K_{3}^{3}=\mu+u(1-\eta)
\beta\left(\frac{1}{\alpha_{1}}+\frac{H_{m}}{\alpha_{0}}\right).
\end{equation}
This completes the proof of Proposition~\ref{proposition2.1.1}.
\end{prooof}\\
\indent Now, consider the following IBVP
\begin{equation}\label{equation2.2}
  \left\lbrace
  \begin{array}{lcr}
    \partial_{t}H-D_{1}\Delta H=f(t,H,I,V)& \mbox{in}& \Omega\times (0,T) \\
   \partial_{t}I-D_{2}\Delta I=g(t,H,I,V) & \mbox{in}&  \Omega\times (0,T) \\
   \partial_{t}V-D_{3}\Delta V=h(t,H,I,V) &\mbox{in}& \Omega\times (0,T)\\
    \frac{\partial H}{\partial\eta}=0; \;\frac{\partial I}{\partial\eta}=0; \; \frac{\partial V}{\partial\eta}=0  &\mbox{on}& \partial\Omega\times [0,T] \\
    H=H_{0}, I=I_{0}, V=V_{0} &\mbox{on}&  \Omega\times\{t=0\}.
  \end{array}
  \right.
\end{equation}
In the following, we will need the following definition and results.
\begin{definition}(Sectorial operator, \cite{Henri1981})
Let A be a linear operator in a Banach space X and suppose A
is closed and densely defined. If there exist real numbers
 $a$, $\omega \in (0, \pi)$,  $M \geq 1$ such that
\begin{equation}\label{sectorial1}
\rho(A)\supset\Sigma=\left\{\lambda_{0}\in \mathbb{C}
  \;:\;\omega\leq arg(\lambda_{0}-a)\leq \pi, \lambda_{0}\neq 0
  \right\}
\end{equation}
and
\begin{equation}\label{sectorial2}
\|R_{\lambda_{0}}(A)\|\leq \frac{M}{|\lambda_{0}-a|} \;\;\mbox{for
all} \;\;\lambda_{0}\in \Sigma,
\end{equation}
then we say that A is sectorial.
\end{definition}
\begin{remark}
The Neumann realization of the Laplacian $A = - \Delta
$, with domain
$$D(A)=\Big\{\omega\in H^{2}(\Omega):
 \frac{\partial\omega}{\partial\eta}=0\Big\}$$
 is a sectorial operator in $L^{2}(\Omega)$.
 But since $C_{0}^{\infty}(\Omega)\subset D(A)$, it is densely
 defined in $L^{2}(\Omega)$. For $\beta\geq0$ large
enough, we define the fractional powers of the Helmholtz operator,
$\mathcal{H^{\beta}}=-\Delta+\beta I$, with domain
$D(\mathcal{H^{\beta}})$ equipped with graph norm
$\|\,.\,\|_{D(\mathcal{H^{\beta}})}=\|\,.\,\|_{2}+\|\mathcal{H^{\beta}}
    . \|_{2}$.
\end{remark}
    We have the following general results.
\begin{lemma}\cite{Adams1975}\label{lem2.1.1}
Let $1\leq p<\infty$. Then $D(\mathcal{H^{\beta}})\subset
C_{0}^{\infty}(\Omega)$ with continuous injection for $\beta >
\frac{n}{2p}$.
\end{lemma}
\begin{lemma}\cite{Henri1981}\label{lem2.1.2}
$D(\mathcal{H}^{\beta})\subset C_{B}^{0}(\Omega)$ with continuous
injection for $\beta
> \frac{n}{4}$.
\end{lemma}
\begin{theorem}\cite{Henri1981}\label{theo2.1.1}
If  A is sectorial, then $-A$ is the infinitesimal
generator of an analytic semigroup, $G(t)$.\\
  If $R_{\lambda_{0}}> a,\; a\in \mathds{R}$ whenever
  $\lambda_{0}\in \sigma$, then for any $t > 0 $,
  $$\|G(t)\|\leq Ce^{-at}\;,\;\;\;\|AG(t)\|
  \leq \frac{C}{t}e^{-at}$$ and
   $$\frac{d}{dt}G(t)=-AG(t)\; ,\;\; t>0.$$
\end{theorem}
\begin{corollary}\label{corol2.1.1}
Let G be the analytic semigroup generated by $-A$. The following
properties hold for the semigroup G and the fractional powers of the
Helmholtz operator $\mathcal{H}^{\beta}$:
  \begin{description}
 \item 1) $G(t)$ : $L^{2}(\Omega)\rightarrow D(\mathcal{H^{\beta}})$
  for all $t>0$,
 \item 2) $\|G(t)\omega\|_\mathcal{H^{\beta}}\leq
  C_{\beta,2}t^{-\beta}\|\omega\|_{2}$ for all $t>0$,
  $\omega \in L^{2}$,
 \item 3) $G(t)\mathcal{H^{\beta}}\omega=\mathcal{H^{\beta}}G(t)\omega$
  for all $t>0$, $\omega \in D(\mathcal{H^{\beta}})$.
  \end{description}
\end{corollary}
\begin{remark}
The following basic hypotheses are assumed to hold :
\begin{description}
  \item[(H1)]$D_{1}>0$,\; $D_{2}>0$ and $D_{3}>0$,
  \item[(H2)] $H_{0}\geq 0$, $I_{0}\geq 0$ and
   $V_{0}\geq 0$ are continuous on $\overline{\Omega}$,
   $H_{0}$, $I_{0}$, $V_{0}$ $\in$ $C_{B}^{0}(\Omega)$,
  \item[(H3)] $f$, $
  g$ and $h$ are continuously differentiable functions from
   $\mathds{\overline{R}}^{4}_{+}$ into $\mathds{R}$ with
   $f(t,0,s,z)\geq 0$, $g(t,r,0,z)\geq 0$ and $h(t,r,s,0)\geq 0$
    for all $t$, $r$, $s$, $z\geq0$.
\end{description}
For $x\in \Omega$, $t\geq0$,  $H$, $I$,
$V\in(C_{B}^{0}(\Omega))^{3}$, define $\mathcal{F},\mathcal{G}$ and
$\mathcal{Q}$ on $\mathds{R_{+}}\times(C_{B}^{0}(\Omega))^{3}$ by :
\begin{center}
$[\mathcal{F}(t,H,I,V)](x)=f(t,H(x),I(x),V(x))$,
$[\mathcal{G}(t,H,I,V)](x)=g(t,H(x),I(x),V(x))$,
$[\mathcal{Q}(t,H,I,V)](x)=h(t,H(x),I(x),V(x))$.
\end{center}
In addition, we let $G_{1}$, $G_{2}$ and $G_{3}$
be the analytical semigroup generated by
 $A_{1}=D_{1}\times\Delta$, $A_{2}=D_{2}\times\Delta$ and $A_{3}=D_{3}\times\Delta$
 respectively.\\
\end{remark}
In the sequel, we will need the following results.
\begin{lemma}\label{lem2.1.3}\cite{Henri1981}
If H, I and V are continuous from $[0, T]$ to $L^{2}(\Omega)$, then
the integrals :
  \begin{center}
   $I_{1}(t)=\int_{0}^{t} G_{1}(t-\tau)\mathcal{F}(\tau,H(\tau),I(\tau),V(\tau))d\tau$,
   $I_{2}(t)=\int_{0}^{t} G_{2}(t-\tau)\mathcal{G}(\tau,H(\tau),I(\tau),V(\tau))d\tau$,
   $I_{3}(t)=\int_{0}^{t} G_{3}(t-\tau)\mathcal{Q}(\tau,H(\tau),I(\tau),V(\tau))d\tau$,
 \end{center}
  exist and  $I_{1}(t)$, $I_{2}(t)$ and $I_{3}(t)$ are continues on
   $[0, T[$ with $I_{1}(t)\in D(A_{1})$, $I_{2}(t)\in D(A_{2})$, $I_{3}(t)\in D(A_{3})$ and
    $I_{1}(t)\rightarrow 0^{+}$ in $L^{2}$ as
     $t \rightarrow 0^{+}$, $I_{2}(t)\rightarrow 0^{+}$
     in $L^{2}$ as $t \rightarrow 0^{+}$ and
     $I_{3}(t)\rightarrow 0^{+}$ in $L^{2}$
  as $t \rightarrow 0^{+}$.
\end{lemma}
\begin{lemma}
If the IBVP (\ref{equation2.2}) has a classical solution, then H, I
and V satisfy the following equalities :
  \begin{eqnarray}
  H(t) &=& G_{1}(t)H_{0}+ \int_{0}^{t} G_{1}(t-\tau)
  \mathcal{F}(\tau,H(\tau),I(\tau),V(\tau))d\tau, \label{equation2.3}\\
  I(t) &=& G_{2}(t)I_{0} + \int_{0}^{t} G_{2}(t-\tau)
  \mathcal{G}(\tau,H(\tau),I(\tau),V(\tau))d\tau,  \label{equation2.4}\\
  V(t) &=& G_{3}(t)V_{0} + \int_{0}^{t} G_{3}(t-\tau)
  \mathcal{Q}(\tau,H(\tau),I(\tau),V(\tau))d\tau. \label{equation2.5}
  \end{eqnarray}
\end{lemma}
\begin{prooof}
Consider the $L^{2}-$valued functions
  $\theta_{j}(\tau)=G_{j}(t-\tau)\omega_{j}(\tau)$,
  $j=1$, $2$, $3$
   with $\omega_{1}=H$, $\omega_{2}=I$ and $\omega_{3}=V$.
   Then $\theta_{j}$ is differentiable since $G_{j}$ is analytic
    and $\omega_{j}$ is differentiable. Then by
    Theorem~\ref{theo2.1.1}, we have
  \begin{eqnarray}\label{equation2.6}
    \frac{d\theta_{1}}{d\tau} &=& \frac{d}{d\tau}\Big[G_{1}(t-\tau)H(\tau)\Big], \nonumber \\
    &=&\frac{d}{d\tau}\Big[G_{1}(t-\tau)\Big]H(\tau) + G_{1}(t-\tau)H^{'}(\tau), \nonumber \\
    &=& -D_{1}\times\Delta\big(G_{1}(t-\tau)\big)H(\tau) + G_{1}(t-\tau)H^{'}(\tau),\nonumber \\
    &=&-D_{1}\times\Delta\big(G_{1}(t-\tau)\big)H(\tau) + G_{1}(t-\tau)\big[D_{1}\times\Delta H(\tau)
       + \mathcal{F}(\tau,H(\tau),I(\tau),V(\tau))\big], \nonumber \\
       &=&-D_{1}\times\Delta\big(G_{1}(t-\tau)\big)H(\tau) +D_{1}\times G_{1}(t-\tau)\times\Delta H(\tau)
       + G_{1}(t-\tau)\mathcal{F}(\tau,H(\tau),I(\tau),V(\tau)), \nonumber
  \end{eqnarray}
  According to corollary~\ref{corol2.1.1} with $\beta=0$, we have
   $$D_{1}\times G_{1}(t-\tau)\Delta H(\tau)=D_{1}\times\Delta
   G_{1}(t-\tau)H(\tau).$$
  \\Therefore
  \begin{eqnarray}\label{equation2.6}
    \frac{d\theta_{1}}{d\tau} &=& -D_{1}\times\Delta G_{1}(t-\tau)H(\tau)
    + D_{1}\times \Delta G_{1}(t-\tau) H(\tau) + G_{1}(t-\tau)\mathcal{F}(\tau,H(\tau),I(\tau),V(\tau)),\nonumber \\
      &=& G_{1}(t-\tau)\mathcal{F}(\tau,H(\tau),I(\tau),V(\tau)).
  \end{eqnarray}
  \begin{eqnarray}\label{equation2.7}
    \frac{d\theta_{2}}{d\tau} &=&\frac{d}{d\tau}\Big[G_{2}(t-\tau)I(\tau)\Big], \nonumber \\
    &=&\frac{d}{d\tau}\Big[G_{2}(t-\tau)\Big]I(\tau) + G_{2}(t-\tau)I^{'}(\tau), \nonumber \\
    &=& -D_{2}\times\Delta\big(G_{2}(t-\tau)\big)I(\tau) + G_{2}(t-\tau)I^{'}(\tau),\nonumber \\
    &=&-D_{2}\times\Delta\big(G_{2}(t-\tau)\big)I(\tau) + G_{2}(t-\tau)\big[D_{2}\times\Delta I(\tau)
       + \mathcal{G}(\tau,H(\tau),I(\tau),V(\tau))\big], \nonumber \\
       &=&-D_{2}\times\Delta\big(G_{2}(t-\tau)\big)I(\tau) +D_{2}\times G_{2}(t-\tau)\times\Delta I(\tau)
       + G_{2}(t-\tau)\mathcal{G}(\tau,H(\tau),I(\tau),V(\tau)), \nonumber \\
    &=& G_{2}(t-\tau)\mathcal{G}(\tau,H(\tau),I(\tau),V(\tau)).
  \end{eqnarray}
  \begin{eqnarray}\label{equation2.8}
    \frac{d\theta_{3}}{d\tau} &=&\frac{d}{d\tau}\Big[G_{3}(t-\tau)V(\tau)\Big], \nonumber \\
      &=&\frac{d}{d\tau}\Big[G_{3}(t-\tau)\Big]V(\tau) + G_{3}(t-\tau)V^{'}(\tau), \nonumber \\
      &=& -D\times\Delta G_{3}(t-\tau)V(\tau)+ G_{3}(t-\tau)V^{ '}(\tau), \nonumber \\
      &=& -D\times\Delta G_{3}(t-\tau)V(\tau)+ G_{3}(t-\tau)\big[D\times\Delta V(\tau)
       + \mathcal{Q}(\tau,H(\tau),I(\tau),V(\tau))\big], \nonumber \\
      &=& -D\times\Delta G_{3}(t-\tau)V(\tau)+ D\times G_{3}(t-\tau)\Delta V(\tau)
      + G_{3}(t-\tau)\mathcal{Q}(\tau,H(\tau),I(\tau),V(\tau)), \nonumber\\
      &=& G_{3}(t-\tau)\mathcal{Q}(\tau,H(\tau),I(\tau),V(\tau)).
  \end{eqnarray}

Integrating equations (\ref{equation2.6}), (\ref{equation2.7})
 and
(\ref{equation2.8}) with respect to time, we obtain equations
(\ref{equation2.3}), (\ref{equation2.4}) and (\ref{equation2.5})
respectively.
\end{prooof}
\begin{remark}
Since in this work, $n=3$, we take $p=2$ so that $\beta>\frac{3}{4}$
and therefore the domain $D(\mathcal{H^{\beta}})$ is continuously
embedded in $C_{0}^{\infty}(\Omega)$ by Lemma~\ref{lem2.1.1}. Now,
let H, I and V be continuous functions from $[0,T]$ to
$D(\mathcal{H^{\beta}})\hookrightarrow C_{B}^{0}(\Omega)$
 satisfying   (\ref{equation2.3}),
 (\ref{equation2.4}) and (\ref{equation2.5}) respectively.
 We can then claim that H, I and V verify system
 (\ref{equation2.2}).
The continuity of  H, I and V implies continuity of
  $t\mapsto \mathcal{F}(t,H(t),I(t),V(t))$,
  $t\mapsto \mathcal{G}(t,H(t),I(t),V(t))$ and
   $t\mapsto \mathcal{Q}(t,H(t),I(t),V(t))$.
   \\\indent
One can then conclude that, the linear Cauchy problem
\begin{equation*}
  \left\lbrace
  \begin{array}{lcr}
    \partial_{t}y_{1}-D_{1}\Delta y_{1}=\mathcal{F}(t,H(t),I(t),V(t)), \\
   \partial_{t}y_{2}-D_{2}\Delta y_{2}=\mathcal{G}(t,H(t),I(t),V(t)), \\
   \partial_{t}y_{3}-D\Delta y_{3}=\mathcal{Q}(t,H(t),I(t),V(t)),\\
   y_{1}(0)=H_{0},y_{2}(0)=I_{0}, y_{3}(0)=V_{0},
  \end{array}
  \right.
\end{equation*}
has a unique solution, with $y_{1}$, $y_{2}$ and $y_{3}$ given by
(\ref{equation2.3}), (\ref{equation2.4}) and (\ref{equation2.5})
respectively.
\end{remark}
Following \cite{Amann1988}, \cite{Amann1990}, \cite{Amann1989},
\cite{Henri1981}, \cite{Hollis1987}, we have the following main
result for the local existence of (\ref{equation1.1}), based on
$L^{2}$-theory.
\begin{proposition}\label{propo3.9}
If hypotheses (H1), (H2) and (H3) are satisfied,
  then the initial value problem (\ref{equation2.2})
  admits a unique solution $(H, I, V)\in (C_{B}^{0}(]0,T],
  D(\mathcal{H}^{\beta})))^{3}$, with $H(0)=H_{0}\in C_{B}^{0}
  (\Omega)$, $I(0)=I_{0}\in C_{B}^{0}(\Omega)$ and $V(0)=V_{0}
  \in C_{B}^{0}(\Omega)$.
\end{proposition}

The proof of this proposition is given in "Appendix A ".

\subsection{Boundedness of the solutions for IBVP
(\ref{equation1.1})}
\begin{proposition}\label{proposition2.2.1}
Let
  $(H,I,V)\in \left( C^{0}\left(\overline{\Omega}\times [0,T)\right)\cap
   C^{2,1}_{B}\left(\Omega\times [0,T)\right) \right)^{3}$\\
   be solutions of (\ref{equation1.1})
   with bounded initial conditions i.e. $0<H_{0}(x)<H_{m}$,
   $0<I_{0}(x)<H_{m}$, $0<V_{0}(x)<V_{m}$
   for all $x\in\overline{\Omega} $,
    and satisfying the boundary condition
    $\frac{\partial H_{0}}{\partial \eta}=0,$ $\frac{\partial I_{0}}{\partial \eta}=0,$ $\frac{\partial V_{0}}{\partial \eta}=0$
    on $\partial \Omega$.
  Then,
\begin{equation*}
\forall (x,t)\in \overline{\Omega}\times [0,T], \;\;H(x,t)\leq
H_{m},\;\;I(x,t)\leq H_{m}\;\;\mbox{and}\;\;V(x,t)\leq V_{m}
\end{equation*}
with
\begin{equation*}
H_{m}=\max\left\{\frac{\lambda}{\delta_{2}},
\max_{x\in \bar{\Omega}} \{H(x,0)+I(x,0)\} \right\}\;\;\mbox{and}\;\;
   V_{m}=\max \left\{ \frac{(1-\varepsilon)kH_{m}}{\mu},
   \max_{x\in \overline{\Omega}}V_{0}(x)\right\}.
\end{equation*}
\end{proposition}

\begin{prooof}
Consider the function $S$ defined  for all $(x,t)\in
\overline{\Omega}\times [0,T]$ by
$$S(x,t)=H(x,t)+I(x,t).$$
 Adding
the first two equations in (\ref{equation1.1}), yields
\begin{equation*}
\frac{\partial S(x,t)}{\partial t} -D_{1}\Delta H(x,t)-D_{2}\Delta I(x,t)= \lambda -dH(x,t) -\alpha
I(x,t).
\end{equation*}
It follows that
  \begin{eqnarray*}
\frac{\partial S(x,t)}{\partial t} -\max\{D_{1},D_{2}\} \Delta\left( H(x,t)+ I(x,t)\right)&\leq& \lambda -\min\{d, \alpha\} \left( H(x,t)+
I(x,t)\right),
\end{eqnarray*}
we have
\begin{equation}\label{equationx.22}
  \left\lbrace
  \begin{array}{lcr}
    \frac{\partial S(x,t)}{\partial t} -\delta_{1} \Delta S(x,t) \leq \lambda -\delta_{2}  S(x,t), \;\; x\in \Omega, \;\; t\in [0,T]\\ \\
   \frac{\partial S(x,t)}{\partial \eta} =0, \;\; x\in \partial\Omega, \;\; t\in [0,T] \\
   S(x,0) =\max\limits_{x\in \overline{\Omega}}S_{0}(x),
  \end{array}
  \right.
\end{equation}
where $S_{0}(x)=\{H(x,0)+I(x,0) \}$, $\delta_{1}=\max\{D_{1},D_{2}\}$ et $\delta_{2} = \min \{d, \alpha \}$. 
 By using the standard parabolic comparison of the scalar
   parabolic equations\cite{MP}, one has
  $$S(x,t) \leq \bar{S}(t),$$
  where
   $\bar{S}(t)=\frac{\lambda}{\delta_{2}} \left(1-e^{-\delta_{2} t}\right)  +\max\limits_{x\in \bar{\Omega}}S_{0}(x)e^{-\delta_{2} t}$
  is the solution of the problem
  \begin{equation}
  \left\lbrace
  \begin{array}{lcr}\label{equationx.23}
    \frac{d \bar{S}(t)}{dt}= \lambda -\delta_{2}\bar{S}(t),\\ \\
   \bar{S}(0) =\max\limits_{x\in \bar{\Omega}}S_{0}(x),
  \end{array}
  \right.
\end{equation}
which dominates system (\ref{equationx.22}).
The general solution of (\ref{equationx.23}) is on the form
 $\bar{S}(t)=k(t)e^{-\delta_{2} t}$. By Lagrange's method, we have
 $k(t)= \frac{\lambda}{\delta_{2}}e^{\delta_{2} t}+c$, $c\in
 \mathbb{R}$.
 Hence $$\bar{S}(t)=\left(\frac{\lambda}{\delta_{2}}e^{\delta_{2} t}+c\right)e^{-\delta_{2} t}.$$
Initial condition yields
  $c=\max\limits_{x\in \overline{\Omega}}S_{0}(x)-\frac{\lambda}{\delta_{2}}$. \\
 Therefore
 $$\bar{S}(t)=\frac{\lambda}{\delta_{2}} \left(1-e^{-\delta_{2} t}\right)  +\max\limits_{x\in \bar{\Omega}}S_{0}(x)e^{-\delta_{2} t}.$$
 Then, it follows that :
 \begin{eqnarray*}
 S(x,t)&\leq & \bar{S}(t)\\
&\leq & \frac{\lambda}{\delta_{2}} \left(1-e^{-\delta_{2} t}\right)  +\max\limits_{x\in \bar{\Omega}}S_{0}(x)e^{-\delta_{2} t}\\
&\leq & \max \left\{\frac{\lambda}{\delta_{2}},\max\limits_{x\in \bar{\Omega}}S_{0}(x)\right\}\left(1-e^{-\delta_{2} t}\right)+ \max \left\{\frac{\lambda}{\delta_{2}},\max\limits_{x\in \bar{\Omega}}S_{0}(x)\right\}e^{-\delta_{2} t}\\
&\leq & \max \left\{\frac{\lambda}{\delta_{2}},\max\limits_{x\in \bar{\Omega}}S_{0}(x)\right\}.
 \end{eqnarray*}
Thus, $$S(x,t)\leq \max\Big\{\frac{\lambda}{\delta_{2}}, \max_{x\in
\bar{\Omega}} \{H(x,0)+I(x,0)\} \Big\}.$$
 Therefore

 $$S(x,t)\leq H_{m}=\max\left\{\frac{\lambda}{\delta_{2}},
\max_{x\in \bar{\Omega}} \{H(x,0)+I(x,0)\} \right\},\forall (x, t)\in \Omega \times [0, T_{max}[,$$
  where $T_{max}$ is the maximal time of existence of the
   solution
   of system (\ref{equation1.1}), this implies that S is bounded.

  Hence $H$ and $I$ are bounded since S is bounded. This prove that
  H and I are bounded.\\\indent
Now, to show that V is bounded, from the third equation of IBVP
 (\ref{equation1.1}), we have
\begin{eqnarray*}
   && \frac{\partial V(x,t)}{\partial t}-D_{3}\Delta V(x,t)
  \leq (1-\varepsilon)kI(x,t)-\mu V(x,t), \;\; x\in \Omega, \;\; t\in [0,T] \\
   && \frac{\partial V(x,t)}{\partial \eta}=0, \;\; x\in \partial\Omega, \;\; t\in [0,T] \\
   && V(x,0)=\max_{x\in \overline{\Omega}}V_{0}(x).
\end{eqnarray*}
It follows from the  previous system, inequality
\begin{equation}\label{equation2.13}
  \left\lbrace
  \begin{array}{lcr}
    \frac{\partial V(x,t)}{\partial t}-D\Delta V(x,t) \leq
    (1-\varepsilon)kH_{m}-\mu V(x,t)\\ \\
   \frac{\partial V(x,t)}{\partial \eta} =0 \\
   V(x,0) =\max\limits_{x\in \bar{\Omega}}V_{0}(x),
  \end{array}
  \right.
\end{equation}
  By using the standard parabolic comparison of the scalar
   parabolic equations \cite{MP}, one has
  $$V(x,t) \leq \overline{V}(t),$$
  where
   $\overline{V}(t)=\frac{(1-\varepsilon)kH_{m}}{\mu}(1-e^{-\mu t}) + \max_{x\in \bar{\Omega}}V_{0}(x)e^{-\mu t}$
  is the solution of the problem
  \begin{equation}\label{equation2.14}
  \left\lbrace
  \begin{array}{lcr}
    \frac{d \overline{V}(t)}{dt}= (1-\varepsilon)kH_{m}-\mu
     \overline{V}(t),\\ \\
   \overline{V}(0) =\max\limits_{x\in \overline{\Omega}}V_{0}(x),
  \end{array}
  \right.
\end{equation}
which dominates system (\ref{equation2.13}).
Indeed, the general solution of
(\ref{equation2.14}) is on the form $\overline{V}(t)=c(t)e^{-\mu
t}$. By the Lagrange's method, we have $c(t)=
\frac{(1-\varepsilon)kH_{m}}{\mu}(e^{\mu t}-1)+ c_{0}$,
$c_{0}\in \mathbb{R}$.\\
thus,
$$\overline{V}(t)=\Big[\frac{(1-\varepsilon)kH_{m}}{\mu}(e^{\mu
t}-1)+c_{0} \Big]e^{-\mu t}.$$ Initial condition yields $\max_{x\in
\bar{\Omega}}V_{0}(x)=\overline{V}(0)=c_{0}$. It follows from that
\begin{equation*}
 \overline{V}(t)=\frac{(1-\varepsilon)kH_{m}}{\mu}(1-e^{-\mu t}) + \max_{x\in \bar{\Omega}}V_{0}(x)e^{-\mu t}.
\end{equation*}
Therefore
 \begin{eqnarray*}
 V(x,t)&\leq & \overline{V}(t)\\
&\leq & \max \Big\{\frac{(1-\varepsilon)kH_{m}}{\mu},\max_{x\in \bar{\Omega}}V_{0}(x)\Big\}(1-e^{-\mu t})+\max \Big\{\frac{(1-\varepsilon)kH_{m}}{\mu},\max_{x\in \bar{\Omega}}V_{0}(x)\Big\}e^{-\mu t}\\
&\leq & \max \Big\{\frac{(1-\varepsilon)kH_{m}}{\mu},\max_{x\in
\bar{\Omega}}V_{0}(x)\Big\}.
 \end{eqnarray*}
  Since
  $$V(x,t) \leq \overline{V}(t) \leq \max
  \Big\{ \frac{(1-\varepsilon)kH_{m}}{\mu},
  \max_{x\in \overline{\Omega}}V_{0}(x)\Big\},
  \;\;\forall (x,t) \in \overline{\Omega}\times [0, T_{max}[;$$
   where $T_{max}$ is the maximal time of existence of the
   solution
   of system (\ref{equation1.1}), this implies that V is bounded.\\
  Thus $H(x,t), I(x,t)$ and $V(x,t)$ are bounded on
   $\overline{\Omega}\times [0, T_{max}[$.
   Therefore, it follows from the standard theory of
     semi-linear parabolic system in \cite{HD1993} that $T_{max}=+\infty$.
     This completes the proof of proposition~\ref{proposition2.2.1}.
\end{prooof}
\subsection{Global existence, uniqueness and positivity
for the IBVP (\ref{equation1.1})} We recast the IBVP
(\ref{equation1.1}) as follows:
\begin{numcases}\strut
\nonumber\frac{\partial w}{\partial t}-\overline{D}\Delta w+
q(w)w=f(w) \;\;\;\mbox{in}\;\; \Omega\times [0, T),\\
 \frac{\partial w_{1}}{\partial \eta} =0,\; \frac{\partial w_{2}}{\partial \eta} =0,\; \frac{\partial w_{3}}{\partial \eta} =0 \;\;\;\mbox{on}\;\;\partial \Omega\times[0, T),
 \\\label{equation2.15}
\nonumber w(x,0) =w_{0}(x)\;\;\;\mbox{in}\;\;\Omega,
\end{numcases}
where $w =(w_{1}, w_{2}, w_{3})^{T}=(H, I, V)^{T}$,
$\overline{D}=diag(D_{1}, D_{2}, D_{3})$, $q(w)=diag\left(q_{1}(w),
q_{2}(w), q_{3}(w)\right)$, $f(w)=\big(f_{1}(w), f_{2}(w),
f_{3}(w)\big)^{T}$, with
\begin{eqnarray*}
q_{1}(w) &=& d+\frac{(1-\eta)\beta w_{3}}{\alpha_{0}+\alpha_{1}w_{1}+\alpha_{2}w_{3}
+\alpha_{3}w_{1}w_{3}},\;\; q_{2}(w) =  (\alpha+\rho),\\ \\
q_{3}(w) &=& \mu + \frac{u(1-\eta)\beta
w_{1}}{\alpha_{0}+\alpha_{1}w_{1}+\alpha_{2}w_{3}
+\alpha_{3}w_{1}w_{3}},\;\;f_{1}(w) = \lambda + \rho w_{2}, \\ \\
 f_{2}(w) &=& \frac{(1-\eta)\beta
w_{1}w_{3}}{\alpha_{0}+\alpha_{1}w_{1}+\alpha_{2}w_{3}
+\alpha_{3}w_{1}w_{3}},\;\;f_{3}(w) = (1-\varepsilon)kw_{2},\\ \\
\end{eqnarray*}
Note that $D_{1}, D_{2}, D_{3}>0$. Denote
  $\mathcal{H}=L^{2}(\Omega)$ and $E=H^{1}(\Omega)$
and define as in \cite{WS} the Hilbert space
 \begin{eqnarray}
 W(0, T, E, E')= \left\{w\in L^{2}\left( (0,T),E \right):\;\;
  \frac{\partial w}{\partial t} \in L^{2}\left( (0,T),E'
  \right) \right\},
\end{eqnarray}
endowed with the norm
\begin{eqnarray}
\|w\|^{2}_{W}=\|w\|^{2}_{L^{2}((0,T),E)} + \Big \|\frac{\partial
w}{\partial t}\Big \|^{2}_{L^{2}((0,T),E')}
\end{eqnarray}
and the following hypothesis for initial conditions:
\begin{eqnarray}\label{equation2.018}
w_{01}\in L^{\infty}(\Omega), w_{02}, w_{03}\in \mathcal{H} \text{
and}\; w_{0i}\geq 0 \text{ for } i\in \{1,2,3\}.
\end{eqnarray}
Here, we apply Theorem~2.7 of \cite{WS}.
So, one approaches the solution by a sequence of solutions of linear
equations. For  $n=0$, $w^{0}$ denotes the solution of
\begin{equation}\label{equation2.18}
  \left\lbrace
  \begin{array}{lcr}
    \frac{\partial w^{0}}{\partial t}-\overline{D}
    \Delta w^{0} = 0 & \mbox{in} & \Omega\times (0, T), \\
  w^{0}(0) =w_{0} & \mbox{in} & \Omega,  \\
   \frac{\partial w^{0}}{\partial \eta} =0. & \mbox{on}
& \partial \Omega
  \end{array}
  \right.
\end{equation}
This equation admits a strong solution and $w^{0}\geq 0$.\\
By induction, $w^{n}$ denotes the solution of
\begin{equation}\label{equation2.19}
  \left\lbrace
  \begin{array}{lcr}
    \frac{\partial w^{n}}{\partial t}-\overline{D}\Delta
    w^{n}+ q(w^{n-1})w^{n}= f(w^{n-1}) & \mbox{in} & \Omega\times (0, T), \\
  w^{n}(0) =w_{0} & \mbox{in} & \Omega,  \\
   \frac{\partial w^{n}}{\partial \eta} =0. & \mbox{on} & \partial
   \Omega.
  \end{array}
  \right.
\end{equation}
Since (\ref{equation2.19}) is a linear equation, $q(w^{n-1})$ and
$f(w^{n-1})$ can replace $a_{0}$ and $f(t)$ of Corollary~2.10 in
\cite{WS}. Suppose that there exists a unique nonnegative solution
$w^{n-1}$. Assuming by induction that $w^{j}\geq 0$ and that by
Proposition~\ref{proposition2.2.1} $w^{j}$   is bounded for $0\leq
j\leq n-1$, one has
 \begin{eqnarray}
0 \leq \frac{ u(1-\eta)\beta
w_{1}^{n-1}}{\alpha_{0}+\alpha_{1}w^{n-1}_{1}
+\alpha_{2}w^{n-1}_{3}+\alpha_{3}w_{1}^{n-1}w^{n-1}_{3}}\leq u(1-\eta)\beta \nonumber\\
\end{eqnarray}
which implies that
 \begin{eqnarray}
 \mu \leq q_{3}(w^{n-1}) \leq \mu + u(1-\eta)\beta.
\end{eqnarray}
Since $w^{j}_{i}$ are bounded , we have
\begin{eqnarray*}
 d \leq q_{1}(w^{n-1}) \leq d + (1-\eta)\beta.
\end{eqnarray*}
In addition, $q_{2}$ is a constant.\\
It then follows that $q_{1}(w^{n-1})$, $q_{2}(w^{n-1})$,
 $q_{3}(w^{n-1}) \in L^{\infty}(\Omega\times(0, T))$.
  We also have $f(w^{n-1})\geq 0$ and $f(w^{n-1})
  \in L^{2}((0, T),
   E')$. Then, by Corollary~2.10 of \cite{WS},
   there exists a unique solution $w^{n}\in W(0, T, E, E')$
   with $w^{n}\geq0$.
Since $ f_{1}(w) = \lambda + \rho w_{2}$, $f_{2}(w) =
\frac{(1-\eta)\beta w_{1}w_{3}}{\alpha_{0}
+\alpha_{1}w_{1}+\alpha_{2}w_{3}+\alpha_{3}w_{1}w_{3}} \leq
(1-\eta)\beta w_{3}$ and $f_{3}(w) = (1-\varepsilon)kw_{2}$, then
$f_{1}(w^{n-1}) =\lambda + \rho w_{2}^{n-1}$,
$f_{2}(w^{n-1})\leq(1-\eta)\beta w_{3}^{n-1} $ and
$f_{3}(w^{n-1})=(1-\varepsilon)kw_{2}^{n-1}$ remain bounded in
$L^{2}(]0,T[, E)$. We deduce that  $w^{n}_{2}$ and $w^{n}_{3}$
remain bounded in $C^{0}([0,T], \mathcal{H})$
 and
$L^{2}((0, T), E)$.
\\\indent
Now, we deduce that the sequence $(w^{n})_{n\geq0}$ (one can extract
a subsequence $(w^{m})_{m\geq0}$) converges weakly to $w_{i}$ in
$L^{2}((0, T), E)$  and weakly star in
 $L^{\infty}((0,T),
\mathcal{H})$ to $w$. Applying Proposition~2.11 in
 \cite{WS} , it holds that for all $n$,
\begin{eqnarray}\label{equation2.023}
 w^{n}(t)&=& G(t)w_{0} + \int_{0}^{t} G(t-s)g^{n}(s)ds,
\end{eqnarray}
where $G(t)$ is the semigroup generated by the unbounded operator
$A= -\overline{D} \Delta$, and
\begin{eqnarray}\label{equation2.024}
g^{n}(s)=-q(w^{n-1}(s))w^{n}(s) + f(w^{n-1}(s)).
\end{eqnarray}
Then, $g^{n} \in L^{2}((0, T), E)$. Since the sequence
$(w^{n})_{n\geq0}$ is bounded in $C^{0}([0,T], \mathcal{H})$, the
sequence $(g^{n})_{n\geq0}$ is bounded in $C^{0}([0,T],
\mathcal{H})$. Now, consider the operator $\mathcal{G}$ from
$C^{0}((0, T), \mathcal{H})$ into $C^{0}((0, T), \mathcal{H})$
defined by
\begin{equation}\label{equation2.025}
 \mathcal{G}(f) =  \int_{0}^{t} G(t-s)f(s)ds.
\end{equation}
Let us prove that  $\mathcal{G}$ is a compact operator. Considering
the triple $(L^{2}(\Omega), H^{1}(\Omega), a)$ with
\begin{eqnarray}
a(w,v)=\sum\limits_{j=1}^{3}\int_{\Omega} \frac{\partial w}{\partial
x_{j}}\frac{\partial v}{\partial x_{j}}dx,
\end{eqnarray}
where $\Omega$ is regular and bounded. As in \cite{WS}, the
unbounded variational operator $A_{\mathcal{H}}$ associated to
 $a$
is a positive symmetric operator with compact resolvent
$R_{\lambda}(A_{\mathcal{H}})$. It admits a sequence
$(\lambda_{k})_{k}$ of positive eigenvalues with
 $\lim\limits_{k\rightarrow
+\infty}\lambda_{k}=+\infty$ and a Hilbert basis $(e_{k})_{k}$ of
$\mathcal{H}$ consisting of eigenvectors of $A_{\mathcal{H}}$. If
$(G(t))_{t>0}$ is the semigroup generated by $-\overline{D}\Delta$, then
for all $w_{0}\in \mathcal{H}$,
\begin{eqnarray}
G(t)w_{0}=\sum_{k=0}^{+\infty}
e^{-t\lambda_{k}}(w_{0},e_{k})e_{k}.
\end{eqnarray}
This proves that the operator is compact for all $t > 0$ since
$$\lim\limits_{k\rightarrow +\infty}e^{-t\lambda_{k}}=0.$$
Setting
\begin{eqnarray}
G_{N}(t)w=\sum_{k=0}^{N} e^{-t\lambda_{k}}(w,e_{k})e_{k},
\end{eqnarray}
one sees that $G_{N}(t)$ is an operator with finite rank which
converges to $G(t)$. The following Theorem is relevant in the
sequel.
\begin{theorem}\label{theo2.4.1}\cite{WS}
Let $t\mapsto G(t)$ be an application from $[0,+\infty)$ into
$\mathcal{L(H)}$. \\
One assumes that there exists a sequence of operators
$(G_{N}(t))_{N\geq0}$ on $\mathcal{H}$ verifying the following
properties:
  \begin{description}
    \item 1) for all $N$ and all $t>0$, $G_{N}(t)$
    has finite rank independent of $t$,
    \item 2)  $t\mapsto G_{N}(t)$ is continuous from
     $[0,+\infty)$ into $\mathcal{L(H)}$ for all $N$,
    \item 3) for $N \rightarrow +\infty$, $G_{N}(t)$
    converges to $G(t)$ in $L^{1}(]0,T[, \mathcal{L(H)})$
    for all $T>0$.
  \end{description}
  Then the operator $\mathcal{G}$ is compact from $C([0,T],
  \mathcal{H})$ to $C([0,T], \mathcal{H})$ for all $T>0$.
\end{theorem}
We are now in the position to prove the global existence, uniqueness
and positivity of the solution to the IBVP (\ref{equation1.1}).
\begin{theorem}\label{theo2.4.2}
 If the initial condition satisfies (\ref{equation2.018}),
 then the IBVP (\ref{equation2.15})
 admits a unique nonnegative
solution $w\in (W(0, T, E, E'))^{3}$.
\end{theorem}

The proof of Theorem~\ref{theo2.4.2} is contained in "appendix B ".

\begin{remark}
It is worth noting that positivity of the solution may be proved by
applying the maximum principle. Moreover, from the above results and
the boundedness of the solution, one has observed that the solution
of IBVP (\ref{equation1.1}) enters the region:
\begin{equation*}
\Sigma =\left\{(H,I,V)\in \Omega^{3}\times\mathds{R}^{3}_{+}:\;
0<H(x,t)\leq H_{m}, 0<I(x,t)\leq H_{m},0<V(x,t)\leq V_{m} \right\},
\end{equation*}
where
\begin{equation*}
 H_{m}=\max\Big\{\frac{\lambda}{\delta_{2}}, \max_{x\in
\overline{\Omega}} \{H(x,0)+I(x,0)\} \Big\}
 \text{ et } V_{m}=\max \Big\{ \frac{(1-\varepsilon)kH_{m}}{\mu},
  \max_{x\in \overline{\Omega}}V_{0}(x)\Big\}.
\end{equation*}
Hence the region $\Sigma$, of biological interest, is
positively-invariant under the flow induced by IBVP
(\ref{equation1.1}).
\end{remark}
\section{Stability analysis of the spatially homogeneous equilibria}\label{sec:3}
\subsection{HCV-spatial homogeneous uninfected equilibrium $E_{0}$}
The spatial homogeneous uninfected equilibrium of the PDE-model
system (\ref{equation1.1}) arises when there is no virus within a
host i.e., V=0. Easy calculations shows that the HCV-spatial
homogeneous uninfected equilibrium for PDE-model system
(\ref{equation1.1}) is given by
$$ E_{0}=(\Lambda, 0, 0) $$
where
$$\Lambda = \frac{\lambda}{d}.$$
\subsection{Basic reproduction number $\mathcal{R}_{0}$}
In order to define the  Basic reproduction number $\mathcal{R}_{0}$
for system (\ref{equation1.1}), we first observe that system
(\ref{equation1.1}) has a spatially homogeneous uninfected
equilibrium $E_{0}$. It should be noted that one of the main tools
in epidemic models is the basic reproduction number
$\mathcal{R}_{0}$ which is an important threshold parameter to
discuss the dynamic behavior of the epidemic model. It quantifies
the infection risk. It measures the expected average number of new
infected hepatocytes generated by a single virion in a completely
healthy hepatocyte. It should be also noted that, while a huge
number of works deals with the threshold dynamics for ODE-models,
very few  studies are devoted to PDE-models. This is eventually due
to the fact that the concept of basic reproduction number has just
recently been extended to PDE-models such as reaction-diffusion and
reaction-convection-diffusion epidemic models with mixed boundary
conditions\cite{Thieme2009,wangzhao2012}. The definition of
$\mathcal{R}_{0}$ in this work follows the approach developed in
\cite{wangzhao2012}.
\\
\par In order to find the basic reproduction number $(\mathcal{R}_{0})$ for the system (\ref{equation1.1}), we obtain the
following linear system at $E_{0}$ for the infected classes:

\begin{numcases}\strut
\nonumber \frac{\partial I}{\partial t}=D_{2}\Delta I-(\alpha+\rho)I+\frac{(1-\eta)\beta\Lambda}{\alpha_{0}+\alpha_{1}\Lambda}V
          \;\; \mbox{in}\;\; \Omega_{T},\\
 \frac{\partial V}{\partial t}=D_{3}\Delta
V+(1-\varepsilon)kI-\mu V-\frac{u(1-\eta)\beta\Lambda}{\alpha_{0}+\alpha_{1}\Lambda}V\;\;
\mbox{in}\;\; \Omega_{T}, \label{syst39}\\
 \nonumber
\frac{\partial I}{\partial\eta}= \frac{\partial V}{\partial\eta}=0  \;\;
\mbox{on} \;\; \partial\Omega\times [0,T].
\end{numcases}
Substituting $I(x,t)=e^{\lambda t}\psi_{2}(x)$ and $V(x,t)=e^{\lambda t}\psi_{3}(x)$ in (\ref{syst39}), we obtain the following cooperative eigenvalue problem:
\begin{numcases}\strut
\nonumber \lambda\psi_{2}(x)=D_{2}\Delta \psi_{2}(x)-(\alpha+\rho)\psi_{2}(x)+\frac{(1-\eta)\beta\Lambda}{\alpha_{0}+\alpha_{1}\Lambda}\psi_{3}(x)
          \;\; \mbox{in}\;\; \Omega,\\
 \lambda\psi_{3}(x)=D_{3}\Delta \psi_{3}(x)+(1-\varepsilon)k\psi_{2}(x)-\mu \psi_{3}(x)-\frac{u(1-\eta)\beta\Lambda}{\alpha_{0}+\alpha_{1}\Lambda}\psi_{3}(x)\;\;
\mbox{in}\;\; \Omega,\; \;\; \;\; \; \label{syst40}\\
 \nonumber
\frac{\partial \psi_{2}(x)}{\partial\eta}= \frac{\partial \psi_{3}(x)}{\partial\eta}=0  \;\;
\mbox{on} \;\; \partial\Omega.
\end{numcases}
As in \cite{wangzhao2012}, let T: $C(\bar{\Omega},\mathbb{R}^{2})\rightarrow C(\bar{\Omega},\mathbb{R}^{2})$ be the solution
semigroup of the following reaction-diffusion system:
\begin{numcases}\strut
\nonumber \frac{\partial I}{\partial t}=D_{2}\Delta I-(\alpha+\rho)I \;\; \mbox{in}\;\; \Omega_{T},\\
 \frac{\partial V}{\partial t}=D_{3}\Delta
V+(1-\varepsilon)kI-\mu V-\frac{u(1-\eta)\beta\Lambda}{\alpha_{0}+\alpha_{1}\Lambda}V\;\;
\mbox{in}\;\; \Omega_{T}, \label{syst41}\\
 \nonumber I(x,0)=\psi_{2}(x), \; V(x,0)=\psi_{3}(x),\; \; \mbox{in}\;\; \Omega_{T}\\
 \nonumber
\frac{\partial I}{\partial\eta}= \frac{\partial V}{\partial\eta}=0  \;\;
\mbox{on} \;\; \partial\Omega.
\end{numcases}
Thus, with initial infection $\Psi (x)=\left(\psi_{2}, \psi_{3}\right)$ the distribution of those infection members
becomes $T(t)\Psi (x)$ as time evolves.
Therefore, the distribution of total new infections is $$ \int_{0}^{\infty} F(x)T(t)\Psi (x)dt, $$
then, we define $$ L(\phi)(x):=\int_{0}^{\infty} F(x)T(t)\Psi (x)dt=F(x) \int_{0}^{\infty} T(t)\Psi (x)dt.$$
$L$ is a positive and continuous operator which maps the initial infection distribution
	 to the distribution of the total infective members produced during the infection
period. Applying the idea of next generation operators \cite{wangzhao2012}, we define the spectral radius
of $L$ as the basic reproduction number $$R_{0}:=\rho (L).$$
The matrices F and V defined as
$$F(x)=
  \begin{pmatrix}
   0 &  \;\; \;\; \frac{(1-\eta)\beta\Lambda}{\alpha_{0}+\alpha_{1}\Lambda} \\ \\
    0 & \;\;\;\; 0  \\ \\
  \end{pmatrix},
\; \;\;
V(x)=
  \begin{pmatrix}
   \alpha+\rho & \;\;\;\;0 \\ \\
    -(1-\varepsilon)k & \;\; \;\; \left[\mu +u\frac{(1-\eta)\beta\Lambda}{\alpha_{0}+\alpha_{1}\Lambda}\right]
  \end{pmatrix}.$$

Then
$$FV^{-1}=\frac{\alpha_{0}+\alpha_{1}\Lambda}{(\alpha+\rho)\left[\mu(\alpha_{0}+\alpha_{1}\Lambda)+ u(1-\eta)\beta\Lambda \right]}
  \begin{pmatrix}
  \frac{(1-\eta)(1-\varepsilon)k\beta\Lambda}{\alpha_{0}+\alpha_{1}\Lambda}  &  \;\; \;\; \frac{(\alpha+\rho)(1-\eta)\beta\Lambda}{\alpha_{0}+\alpha_{1}\Lambda} \\ \\
    0 & \;\;\;\; 0  \\ \\
  \end{pmatrix}.$$
By \cite{wangzhao2012} (theorem~3.4), one has
\begin{equation}\label{basicnumber}
\mathcal{R}_{0}=\frac{(1-\eta)(1-\varepsilon)k\beta\Lambda}
 {(\alpha+\rho)\left[\mu(\alpha_{0}+\alpha_{1}\Lambda)+ u(1-\eta)\beta\Lambda \right]}.
\end{equation}
\subsection{Existence and uniqueness of HCV-spatial homogeneous infected equilibrium $E^{*}$}
In this section, we address the existence and uniqueness of infected
spatial homogeneous equilibrium (\ref{equation1.1}). The latter
denoted as $E^{*}=(H^{*}, I^{*}, V^{*})$ with $H^{*}\neq 0$,
$I^{*}\neq 0$ et $V^{*}\neq 0$ satisfying the following algebraic
system :
\begin{equation}\label{equation3.2}
  \left\lbrace
  \begin{aligned}
    \lambda-dH^{*}-(1-\eta)L(H^{*}, I^{*}, V^{*})V^{*} + \rho I^{*} &=&0,\\
   (1-\eta)L(H^{*}, I^{*}, V^{*})V^{*} -(\alpha+\rho)I^{*} &=&0, \\
  (1-\varepsilon)kI^{*} - \mu V^{*} - u(1-\eta)L(H^{*}, I^{*}, V^{*})V^{*} &=&0,
  \end{aligned}
  \right.
\end{equation}
where
$$L(H, I, V)=\frac{\beta
H}{\alpha_{0}+\alpha_{1}H+\alpha_{2}V+\alpha_{3}HV}.$$ Adding the
first and second equation of (\ref{equation3.2}), we have
$$\lambda-dH^{*}-\alpha I^{*}=0$$ which yields
\begin{equation}\label{ietoile}
I^{*}=\frac{\lambda-dH^{*}}{\alpha}.
\end{equation}
As far as, using the second and third equation of
(\ref{equation3.2}), one has
\begin{eqnarray*}
 -u(\alpha+\rho)I^{*}+(1-\varepsilon)kI^{*}- \mu V^{*}=0 ,
\end{eqnarray*}
i.e.,
\begin{eqnarray*}
 V^{*}=\frac{(1-\varepsilon)k-u(\alpha+\rho)}{\mu}I^{*},
 \end{eqnarray*}
hence
\begin{eqnarray}\label{equation3.3}
V^{*}=\frac{(1-\varepsilon)k-u(\alpha+\rho)}{\mu}
\frac{\lambda-dH^{*}}{\alpha}
 \end{eqnarray}
 according to (\ref{ietoile}).
 The substitution of (\ref{equation3.3}) in (\ref{equation3.2}) yields :
 \begin{eqnarray*}
 (1-\eta)L\left(H^{*},\frac{\lambda-dH^{*}}{\alpha},\frac{(1-\varepsilon)k-u(\alpha+\rho)}{\mu}
 \frac{\lambda-dH^{*}}{\alpha}\right)\frac{(1-\varepsilon)k-u(\alpha+\rho)}{\mu}I^{*}-(\alpha+\rho)I^{*}=0.
 \end{eqnarray*}
 Thus, we have
  \begin{eqnarray*}
  (1-\eta)\Big((1-\varepsilon)k-u(\alpha+\rho)\Big)L\Big(H^{*},
  \frac{\lambda-dH^{*}}{\alpha},\frac{(1-\varepsilon)k-u(\alpha+\rho)}
  {\mu}\frac{\lambda-dH^{*}}{\alpha}\Big) = (\alpha+\rho)\mu
 \end{eqnarray*}
since  $I^{*}\neq0$. Furthermore,
 $I^{*}\geq0$, gives $\frac{\lambda-dH^{*}}{\alpha}\geq0$. Thus $H^{*}\leq \frac{\lambda}{d}$.
 Hence there is not a biological equilibrium when $H^{*}>\frac{\lambda}{d}$.\\
Let us consider the function  $\psi$ defined on $\Big[0,
\frac{\lambda}{d}\Big]$ by :
 \begin{eqnarray*}
\psi (x) &=&(1-\eta)\gamma L\Big(x,\frac{\lambda-d
x}{\alpha},\frac{\gamma (\lambda-d x)}{\mu \alpha}\Big)
-(\alpha+\rho)\mu,
 \end{eqnarray*}
where
$$\gamma =(1-\varepsilon)k-u(\alpha+\rho).$$
We have  $$\psi (0)=-(\alpha+\rho)\mu <0$$ and
 \begin{eqnarray*}
 \psi \left(\frac{\lambda}{d}\right) &=& (1-\eta)\gamma L\left(\frac{\lambda}{d},0,0 \right) -(\alpha+\rho)\mu, \\
 &=& (1-\eta)\left[(1-\varepsilon)k-u(\alpha+\rho) \right] \frac{\beta \Lambda}{\alpha_{0}+\alpha_{1}\Lambda} -(\alpha+\rho)\mu, \\
 &=& \frac{(1-\eta)(1-\varepsilon)k\beta \Lambda}{\alpha_{0}+\alpha_{1}\Lambda} -\frac{u(1-\eta)(\alpha+\rho)\beta \Lambda}{\alpha_{0}+\alpha_{1}\Lambda}-(\alpha+\rho)\mu, \\
 &=& \frac{1}{\alpha_{0}+\alpha_{1}\Lambda} \left[ (1-\eta)(1-\varepsilon)k\beta \Lambda - \mu(\alpha+\rho)(\alpha_{0}+\alpha_{1}\Lambda) -u(1-\eta)(\alpha+\rho)\beta \Lambda \right],\\
 &=& \frac{1}{\alpha_{0}+\alpha_{1}\Lambda} \Big[ (1-\eta)(1-\varepsilon)k\beta \Lambda - (\alpha+\rho)\left[\mu(\alpha_{0}+\alpha_{1}\Lambda) +u(1-\eta)\beta \Lambda\right] \Big],\\
 &=& \frac{(\alpha+\rho)\left[\mu(\alpha_{0}+\alpha_{1}\Lambda) +u(1-\eta)\beta \Lambda\right]}{\alpha_{0}+\alpha_{1}\Lambda}( \mathcal{R}_{0} - 1),
 \end{eqnarray*}

 It follows that
\begin{eqnarray*}
 \psi \Big(\frac{\lambda}{d}\Big) &=& \frac{(\alpha+\rho)\left[\mu(\alpha_{0}+\alpha_{1}\Lambda) +u(1-\eta)\beta \Lambda\right]}{\alpha_{0}+\alpha_{1}\Lambda}( \mathcal{R}_{0} - 1)>0 \;\;
 \text{if and only if} \;\; \mathcal{R}_{0}>1.
 \end{eqnarray*}
\par Moreover, letting
 $y=\dfrac{\lambda-d.x}{\alpha}$ and $z=\dfrac{\gamma (\lambda-d.x)}{\mu \alpha}$,
 we have
 \begin{eqnarray*}
 \psi^{'} (x) &=& (1-\eta)\gamma.\frac{d}{dx}
 \Big[L\left(x,\frac{\lambda-dx}{\alpha},\frac{\gamma (\lambda-dx)}{\mu \alpha}\right)
 -(\alpha+\rho)\mu\Big],\\
 &=& (1-\eta)\gamma \left( \frac{\partial L}{\partial x}-\frac{d}{\alpha}\frac{\partial L}{\partial y}
 -\frac{\gamma d}{\mu \alpha}\frac{\partial  L}{\partial z} \right),\\
 &=& (1-\eta)\gamma \left(\frac{\partial L}{\partial x} - \frac{d}{\alpha}
  \frac{\partial L}{\partial x}\frac{\partial x}{\partial y} -
  \frac{\gamma d}{\mu \alpha}\frac{\partial  L}{\partial x}\frac{\partial x}{\partial z} \right),\\
 &=& (1-\eta)\gamma \left(\frac{\partial L}{\partial x} - \frac{d}{\alpha}\frac{\partial L}{\partial x}
 (- \frac{\alpha}{d}) - \frac{\gamma d}{\mu \alpha}\frac{\partial  L}{\partial x}(- \frac{\mu \alpha}
 {\gamma d}) \right),\\
 &=& 3(1-\eta)\gamma\frac{\partial L}{\partial x},\\
 &=& 3(1-\eta)\gamma\frac{\beta \alpha_{0}+\beta \alpha_{2}V}{\left(\alpha_{0}+(\alpha_{1}+\alpha_{3}V)x
 +\alpha_{2}V \right)^{2}}>0 \;\; \text{ if and only if }\;\; \gamma >0.
 \end{eqnarray*}
Therefore, if $\mathcal{R}_{0}>1$ there exists a unique spatially
homogeneous infected equilibrium $E^{*}=(H^{*}, I^{*}, V^{*})$ with
$H^{*}\in \left(0,\frac{\lambda}{d}\right)$, $I^{*}>0$ and
$V^{*}>0$.
\par
The previous investigations can be summarized in the following
theorem :
\begin{theorem}.\label{theo2.5.1} \\
 1) if $\mathcal{R}_{0}\leq 1$, then the PDE-system (\ref{equation1.1}) admits a unique spatially homogeneous
 uninfected equilibrium
  $E_{0}=\left(\frac{\lambda}{d},0,0\right)$.\\
 2) If  $\mathcal{R}_{0}> 1$, then the PDE system
 (\ref{equation1.1}) admits a unique spatially homogeneous infected equilibrium
 $E^{*}=(H^{*}, I^{*}, V^{*})$ with $H^{*}\in \left(0,\frac{\lambda}{d}\right)$, $I^{*}>0$ and $V^{*}>0$.
\end{theorem}
\begin{remark}
Due to the spacial dependence of the state variables, spatially-inhomogeneous steady states can exist.
\par Indeed, any spatially-inhomogeneous equilibrium point $E=(H,I,V)$ of the model (\ref{equation1.1}) subject to the homogeneous Neumann boundary condition must solve the following system.
\begin{numcases}\strut
\nonumber D_{1}\Delta H+\lambda-dH-\frac{(1-\eta)\beta HV}
{\alpha_{0}+\alpha_{1}H+\alpha_{2}V+\alpha_{3}HV}+\rho I=0,\\
\nonumber D_{2}\Delta I+\frac{(1-\eta)\beta
HV}{\alpha_{0}
          +\alpha_{1}H+\alpha_{2}V+\alpha_{3}HV}-(\alpha+\rho)I=0,\\
 D_{3}\Delta
V+(1-\varepsilon)kI-\mu V-\frac{u(1-\eta)\beta
HV}{\alpha_{0}+\alpha_{1}H+\alpha_{2}V+\alpha_{3}HV}=0,\\
 \nonumber
\frac{\partial H}{\partial\eta}=\frac{\partial I}{\partial\eta}=\frac{\partial V}{\partial\eta}=0. \\
\nonumber
\end{numcases}
After the existence given by theorem~\ref{theo2.5.1}, investigation of the local stability of such spatially-inhomogeneous equilibria will be the concern of the two forthcoming sections.
\end{remark}
\subsection{Local stability of HCV-uninfected equilibrium}
 The objective of  this section is to discuss the local stability of the
spatially homogeneous uninfected equilibrium for the PDE system
(\ref{equation1.1}). We address local stability by analyzing the
characteristic equation.
\begin{theorem}\label{theo3.2.1}
The spatially homogeneous uninfected equilibrium $E_{0}$ of
PDE-model system (\ref{equation1.1}) is locally asymptotically
stable if $\mathcal{R}_{0}\leq 1$ and it is unstable if
$\mathcal{R}_{0} > 1$.
\end{theorem}
\begin{prooof}
Let $\{\mu_{l}, \varphi_{l}\}$ be an eigenpair of the Laplace
operator $-\Delta$ on $\Omega$ with the homogeneous Neumann boundary
condition where $0=\mu_{1}<\mu_{2}<\mu_{3}< \cdot\cdot\cdot$. Let
$E_{\mu_{l}}$ be the eigenspace corresponding to  $\mu_{l}$ in
 $C^{1}(\Omega)$ and $\{ \varphi_{lj}, j=1,2,\cdot\cdot\cdot,\mbox{dim}
E_{\mu_{l}} \}$
 be an orthogonal basis of  $E_{\mu_{l}}$. Let
 $\mathbb{X}=(C^{1}(\Omega))^{3}$ and $\mathbb{X}_{lj}=\{ \varphi_{lj} c,\;/
 \; c\in \mathbb{R}^{3} \}$.\\
Consider the following direct sum
$$ \mathbb{X}=\bigoplus\limits_{l=1}^{\infty}\mathbb{X}_{l}
\;\; \mbox{with}\;\; \mathbb{X}_{l}=\bigoplus\limits_{j=1}^{dim
E_{\mu_{l}}}\mathbb{X}_{lj},$$ where $\mathbb{X}_{lj}$ is the
eigenspace corresponding to $\mu_{l}$. Linearizing
(\ref{equation1.1}) at the spatially homogeneous uninfected
equilibrium $E_{0}$ we obtain the following linearized system :
\begin{numcases}\strut
\nonumber \frac{\partial w_{1}}{\partial t}=D_{1}\Delta w_{1} - dw_{1} + \rho w_{2}
-\frac{(1-\eta)\beta\Lambda}{\alpha_{0}
    +\alpha_{1}\Lambda} w_{3},\\
         \frac{\partial w_{2}}{\partial t}=D_{2}\Delta w_{2} - (\alpha+\rho)w_{2} + \frac{(1-\eta)\beta\Lambda}{\alpha_{0}
   +\alpha_{1}\Lambda} w_{3},\label{equation3.4}\\
\nonumber \frac{\partial w_{3}}{\partial t}=D_{3}\Delta w_{3} +
(1-\varepsilon)k w_{2}- \Big[\mu
   +
   u\frac{(1-\eta)\beta\Lambda}{\alpha_{0}+\alpha_{1}\Lambda}\Big]w_{3},
\end{numcases}
where $W=(w_{1},w_{2},w_{3})^{T}=(H,I,V)^{T}$.\\
From the previous system (\ref{equation3.4}), we obtain
\begin{equation*}
W_{t}=\mathcal{L}W=\overline{D}\Delta W + \mathcal{K}(E_{0})W
\end{equation*}
where $\overline{D}=diag(D_{1},D_{2},D_{3})$ and

\begin{equation}\label{equation3.55}
\mathcal{K}(E_{0})W=
\begin{pmatrix}
 -dw_{1} + \rho w_{2} -\frac{(1-\eta)\beta\Lambda}{\alpha_{0}+\alpha_{1}\Lambda} w_{3} \\
  -(\alpha+\rho)w_{2} + \frac{(1-\eta)\beta\Lambda}{\alpha_{0}+\alpha_{1}\Lambda} w_{3} \\
  (1-\varepsilon)kw_{2}- \left(\mu + u\frac{(1-\eta)\beta\Lambda}
   {\alpha_{0}+\alpha_{1}\Lambda}\right)w_{3} \\
\end{pmatrix}.
\end{equation}
For each $l\geq 1$, $\mathbb{X}_{l}$ is invariant under the operator
$\mathcal{L}$, and $\tilde{\lambda}$ is an eigenvalue of
$\mathcal{L}$ if and only if it is an eigenvalue of the matrix
$-\mu_{l} \overline{D} + \mathcal{K}(E_{0})$ for some $l\geq 1$, in
which case, there is an eigenvector in $\mathbb{X}_{l}$. So, one has
{\small
\begin{equation*}\label{equation3.6}
det\left(-\mu_{l} \overline{D} +
\mathcal{K}(E_{0})-\tilde{\lambda}Id\right)=
  \left|
  \begin{matrix}
    -\left(\mu_{l}D_{1}+d+\tilde{\lambda}\right) & \rho &       -\frac{(1-\eta)\beta\Lambda}{\alpha_{0}+\alpha_{1}\Lambda}  \\ \\
   0   & -\left(\mu_{l}D_{2}+(\alpha+\rho)+\tilde{\lambda}\right) & \frac{(1-\eta)\beta\Lambda}{\alpha_{0}+\alpha_{1}\Lambda} \\ \\
   0   & (1-\varepsilon)k  &  - \left(\mu_{l}D_{3}+\mu +u\frac{(1-\eta)\beta\Lambda}{\alpha_{0}
   +\alpha_{1}\Lambda}\right)
   -\tilde{\lambda}
  \end{matrix}
  \right|
\end{equation*}
}
The characteristic equation of $-\mu_{j}\overline{D}+\mathcal{K}(E_{0})\; $\;is

{\normalsize
\begin{equation}\label{equation3.7}
-(\mu_{l}D_{1}+d+\tilde{\lambda})\left[ \left[\mu_{l}D_{2}+ (\alpha+\rho)+\tilde{\lambda}
\right]\left[ \left(\mu_{l}D_{3}+\mu
+u\frac{(1-\eta)\beta\Lambda}{\alpha_{0}+\alpha_{1}\Lambda}\right)+\tilde{\lambda}
\right]-
\frac{(1-\varepsilon)(1-\eta)k\beta\Lambda}{\alpha_{0}+\alpha_{1}\Lambda}\right]=0,
\end{equation}
} from (\ref{equation3.7}), we get
$$\tilde{\lambda}_{0}=-\mu_{l}D_{1}-d<0,\; and$$
another characteristic eigenvalues are the roots of the following equation :
\begin{eqnarray}\label{equation3.8}
\tilde{\lambda}^{2}+B\tilde{\lambda} +
(\mu_{l}D_{2}+\alpha+\rho)\left(\mu_{l}D_{3}+\mu
+u\frac{(1-\eta)\beta\Lambda}{\alpha_{0}+\alpha_{1}\Lambda}\right) -
\frac{(1-\varepsilon)(1-\eta)k\beta\Lambda}{\alpha_{0}+\alpha_{1}\Lambda}=0,\;\;\;\;\;\;
\end{eqnarray}
where $$B=\left[ \left(\mu_{l}D_{3}+\mu
+u\frac{(1-\eta)\beta\Lambda}{\alpha_{0}+\alpha_{1}\Lambda}\right)+\mu_{l}D_{2}+
(\alpha+\rho)\right].$$

 $$Let \; \; C=(\mu_{l}D_{2}+\alpha+\rho)\left(\mu_{l}D_{3}+\mu
+u\frac{(1-\eta)\beta\Lambda}{\alpha_{0}+\alpha_{1}\Lambda}\right) -
\frac{(1-\varepsilon)(1-\eta)k\beta\Lambda}{\alpha_{0}+\alpha_{1}\Lambda}.$$
One has,
 \begin{eqnarray*}
 C &=& \mu_{l}D_{2}\left(\mu_{l}D_{3}+\mu
+u\frac{(1-\eta)\beta\Lambda}{\alpha_{0}+\alpha_{1}\Lambda}\right) +(\alpha+\rho)\mu_{l}D_{3}+ (\alpha+\rho)\left(\mu
+u\frac{(1-\eta)\beta\Lambda}{\alpha_{0}+\alpha_{1}\Lambda}\right)\\
 &-&\frac{(1-\varepsilon)(1-\eta)k\beta\Lambda}{\alpha_{0}+\alpha_{1}\Lambda},\\
&=& \mu_{l}D_{2}\left(\mu_{l}D_{3}+\mu
+u\frac{(1-\eta)\beta\Lambda}{\alpha_{0}+\alpha_{1}\Lambda}\right) +(\alpha+\rho)\mu_{l}D_{3}\\
&+& \frac{1}{\alpha_{0}+\alpha_{1}\Lambda}\Big[(\alpha+\rho)\left[\mu(\alpha_{0}+\alpha_{1}\Lambda)
+u(1-\eta)\beta\Lambda\right] -
(1-\varepsilon)(1-\eta)k\beta\Lambda\Big],\\
&=& \mu_{l}D_{2}\left(\mu_{l}D_{3}+\mu+u\frac{(1-\eta)\beta\Lambda}{\alpha_{0}+\alpha_{1}\Lambda}\right) +(\alpha+\rho)\mu_{l}D_{3}\\
&+&\frac{(\alpha+\rho)\left[\mu(\alpha_{0}+\alpha_{1}\Lambda)+ u(1-\eta)\beta\Lambda\right]}{\alpha_{0}+\alpha_{1}\Lambda}(1 -\mathcal{R}_{0}).
 \end{eqnarray*}

Since $B > 0$, if $\mathcal{R}_{0}\leq 1$ then $C$ is also positive. Hence by vertue of the Routh-Hurwitz
criterion, equation (\ref{equation3.8}) does not admit solution with
positive real part. Thus none characteristic eigenvalue have
positive real part. Therefore if $\mathcal{R}_{0}\leq 1$, the
spatially homogeneous uninfected equilibrium
$E_{0}=\Big(\frac{\lambda}{d},0,0\Big)$ of (\ref{equation1.1}) is
locally asymptotically stable.
\par Otherwise if $\mathcal{R}_{0}> 1$, then for $l=1$, (in this case $\mu_{1}=0$) one has,
$$C=\frac{(\alpha+\rho)\left[\mu(\alpha_{0}+\alpha_{1}\Lambda)
+u(1-\eta)\beta\Lambda\right]}{\alpha_{0}+\alpha_{1}\Lambda}(1 -\mathcal{R}_{0})<0.$$
Hence there is a complexe root of equation (\ref{equation3.8}) with
positive real part in the spectrum of $\mathcal{K}$ according to
Routh-Hurwitz criterion. Therefore the uninfected equilibrium
$E_{0}=\Big(\frac{\lambda}{d},0,0\Big)$ of
 (\ref{equation1.1}) is  unstable. This completes the proof of
 theorem~\ref{theo3.2.1}.
\end{prooof}
\subsection{Global stability of HCV-Uninfected equilibrium}
The objective of  this section is to discuss the global stability of
the spatially homogeneous uninfected equilibrium for the PDE system
(\ref{equation1.1}). We address global stability by using the method
of construction of Lyapunov functionals. These Lyapunov functional
is obtained from those  differential equations by applying the
method of Hattaf and Yousfi presented in \cite{hattafyousfi2013}.
For this purpose, we start by letting
\begin{equation*}
  \tau_{0}=\frac{(1-\varepsilon)k(1-\eta)\beta \Lambda}
  {\mu\alpha_{0}(\alpha+\rho)}.
\end{equation*}
Then, it is easy to see that
\begin{equation*}
 \frac{(1-\varepsilon)(1-\eta)k\beta\Lambda}
 {(\alpha+\rho)\left[\mu(\alpha_{0}+\alpha_{1}\Lambda)
 +u(1-\eta)\beta\Lambda\right]} \leq \frac{(1-\varepsilon)
 k(1-\eta)\beta \Lambda}{\mu\alpha_{0}(\alpha+\rho)},
\end{equation*}
i.e.,
 $$\mathcal{R}_{0}\leq\tau_{0}.$$
 We state the following result on global stability at $E_{0}$ as
 follows :
\begin{theorem}\label{theogloE0}
The spatially homogeneous uninfected equilibrium $E_{0}$ of
PDE-model system (\ref{equation1.1}) is globally asymptotically
stable in the positively-invariant region $\Sigma$ if
$\tau_{0}<1$.
\end{theorem}
\begin{prooof}
Let us consider the following function
  \begin{equation*}
    G_{1}(t)=\frac{(1-\varepsilon)k}{\alpha+\rho}I(t) + V(t).
  \end{equation*}
Then, the differention of $G_{1}$ with respect to $t$ gives
\begin{equation*}
    \frac{dG_{1}}{dt} = \left(\frac{(1-\varepsilon)k - u(\alpha + \rho)}{\mu (\alpha + \rho)
    (\alpha_{0}+\alpha_{1}H+\alpha_{2}V+\alpha_{3}HV)}(1-\eta)\beta H-1\right) \mu V.
  \end{equation*}
Since $H \leq \frac{\lambda}{d}=\Lambda$ in the positively-invariant
region $\Sigma$, one has
\begin{eqnarray*}
  \frac{dG_{1}}{dt}  &\leq& \Bigg[\frac{\big[(1-\varepsilon)k - u(\alpha
  + \rho)\big](1-\eta)\beta \Lambda}{\mu\alpha_{0}(\alpha + \rho)}-1\Bigg] \mu V,\\
  &\leq& \Big[\frac{(1-\varepsilon)(1-\eta)k\beta \Lambda}{\mu\alpha_{0}(\alpha + \rho)}-1\Big] \mu V,\\
    &\leq& (\tau_{0}-1)\mu V.
  \end{eqnarray*}
Now, we define the Lyapunov function as follows
  $$L_{1}=\int_{\Omega} G_{1}dx. $$
The computation of the time derivative of $L_{1}$ along the positive
solutions of the PDE-model system (\ref{equation1.1}) yields
\begin{eqnarray*}
  \frac{dL_{1}}{dt} &=& \frac{d}{dt}\Big[ \int_{\Omega} G_{1}dx \Big],\\
   &=& \int_{\Omega} \frac{dG_{1}}{dt}dx, \\
   &\leq& \int_{\Omega} \Big[ (\tau_{0}-1)\mu V \Big] dx.
\end{eqnarray*}
It is clear that the condition $\tau_{0}\leq 1$ gives
$\frac{dL_{1}}{dt}\leq 0$ for all H, I, V $>$ 0. We note that the
solutions of system (\ref{equation1.1}) are limited by $\Upsilon$,
the greatest invariant subset of E=$\Big\{ (H, I, V)\in \Sigma |
\frac{dL_{1}}{dt}=0 \Big\}$. We realize that $\frac{dL_{1}}{dt}=0$
if and only if V = 0 and I = 0. Each element of $\Upsilon$ satisfies
V = 0 and consequently I=0. By Lyapunov-LaSalle invariance principle
\cite{khalil1}, $E_{0}$ is globally asymptotically stable if
$\tau_{0}< 1$. So, we obtain a sufficient condition $\mathcal{R}_{0}
\leq \tau_{0}$ which ensures that the HCV spatially homogeneous
equilibrium $E_{0}$ of PDE-model system (\ref{equation1.1}) is
globally asymptotically stable if $\tau_{0}< 1$. This completes the
proof of theorem~\ref{theogloE0}.
\end{prooof}
\subsection{Local stability of HCV spatially
homogeneous infected equilibrium}
Let us study the local stability of the unique infected spatially
homogeneous equilibrium $E^{*}$ of our PDE-model system. Consider
the Laplace operator $-\Delta$ and let $0=\mu_{1}<\mu_{2}<\mu_{3}<
\cdot\cdot\cdot$ be its eigenvalues on $\Omega$ with the homogeneous
Neumann boundary condition, and
 $E_{\mu_{l}}$ be the eigenspace corresponding to  $\mu_{l}$ in
 $C^{1}(\Omega)$. Let also
 $\mathbb{X}=(C^{1}(\Omega))^{3}$,
  $\{\varphi_{lj}, j=1,2,\cdot\cdot\cdot,dim
E_{\mu_{l}} \}$
 be an orthogonal basis of  $E_{\mu_{l}}$
 and $\mathbb{X}_{lj}=\{ \varphi_{lj} c\;/\; c\in \mathbb{R}^{3} \}$.\\
Then, $$ \mathbb{X}=\bigoplus\limits_{l=1}^{\infty}\mathbb{X}_{l}
\;\; \mbox{with}\;\; \mathbb{X}_{l}=\bigoplus\limits_{j=1}^{dim
E_{\mu_{l}}}\mathbb{X}_{lj}.$$
 Now, let set
$w_{1}=H$, $w_{2}=I$, $w_{3}=V$. Further we use the vecteur notation
$W=(w_{1},w_{2},w_{3})^{T}=(H,I,V)^{T}$, and
$\overline{D}=diag(D_{1},D_{2},D_{3})$. Then the linearization of  the PDE
system  at $E^{*}$ is of the form
$$W_{t}=\mathcal{L}W=\overline{D}\Delta W + \mathcal{K}(E^{*})W,$$
where
\begin{equation}\label{equation3.12}
\mathcal{K}(E^{*})W=
\begin{pmatrix}
  -\left( d+ A \right)w_{1} + \rho w_{2} -B w_{3} \\
  A w_{1}-(\alpha+\rho)w_{2} + B w_{3} \\
  -u A w_{1}+(1-\varepsilon)kw_{2}- \left(\mu + u B\right)w_{3}\\
\end{pmatrix},
\end{equation}
with
\begin{equation*}
  A=\frac{(1-\eta)(\alpha_{0}+\alpha_{2}V^{*})\beta V^{*}}{(\alpha_{0}+\alpha_{1}H^{*}+\alpha_{2}V^{*}+\alpha_{3}H^{*}V^{*})^{2}}
\end{equation*}
and
\begin{equation*}
  B=\frac{(1-\eta)(\alpha_{0}+\alpha_{1}H^{*})\beta H^{*}}{(\alpha_{0}
  +\alpha_{1}H^{*}+\alpha_{2}V^{*}+\alpha_{3}H^{*}V^{*})^{2}}.
\end{equation*}
For each  $l\geq 1$, $\mathbb{X}_{l}$ is invariant under the
operator  $\mathcal{L}$, and $\tilde{\lambda}$ is an eigenvalue
$\mathcal{L}$ if and only if it is an eigenvalue of the matrix
$-\mu_{l} \overline{D} + \mathcal{K}(E^{*})$ for some $l\geq 1$, in
which case, there is an eigenvector in $\mathbb{X}_{l}$. Therefore
we get :
\begin{equation*}\label{equation3.13}
det\Big(-\mu_{l} \overline{D} +
\mathcal{K}(E^{*})-\tilde{\lambda}Id\Big)=
  \left|
  \begin{matrix}
    -( \mu_{l}D_{1}+d+ A)-\tilde{\lambda} & \rho &       -B \\
   A   & -(\mu_{l}D_{2}+\alpha+\rho)-\tilde{\lambda} & B \\
   -u A   & (1-\varepsilon)k  &  - \left(\mu_{l}D_{3} + \mu +u B\right)-\tilde{\lambda}
  \end{matrix}
  \right|.
\end{equation*}
The characteristic equation of  $-\mu_{l} \overline{D} +
\mathcal{K}(E^{*})$ is on the form
\begin{equation}\label{equation3.14}
  \tilde{\lambda}^{3} + a_{2} \tilde{\lambda}^{2} + a_{1}\tilde{\lambda}+a_{0}=0
\end{equation}
where
\begin{equation*}
a_{2} = (\mu_{l}D_{1}+d+A+\mu_{l}D_{2}+\alpha +\rho +\mu_{l}D_{3} + \mu + u B)>0,
\end{equation*}
\begin{equation*}
a_{1} =(\mu_{l}D_{1}+d+A)(\mu_{l}D_{2}+\alpha +\rho +\mu_{l}D_{3} + \mu + u B)+(\mu_{l}D_{2}+\alpha +\rho)(\mu_{l}D_{3} + \mu + u B)-(1-\varepsilon)kB,
\end{equation*}
\begin{equation*}
a_{0} =(\mu_{l}D_{1}+d+A)(\mu_{l}D_{2}+\alpha +\rho)(\mu_{l}D_{3} + \mu + u B)- (\mu_{l}D_{1}+d+A)(1-\varepsilon)kB.
\end{equation*}
If  $a_{1}>0$ and $a_{1} a_{2}>a_{0}$ from the above investigations,
it then follows from Routh-Hurwitz criterion that all roots of
(\ref{equation3.14}) have negative real parts and therefore we have
the following result.
\begin{theorem}\label{theo3.4.1}  If  $a_{1}>0$ and $a_{1} a_{2}>a_{0}$,
 then the spatially homogeneous
infected equilibrium  $E^{*}=(H^{*}, I^{*}, V^{*})$ of the PDE-model
system (\ref{equation1.1}) is locally asymptotically stable when it
exists.
\end{theorem}
\subsection{Global stability of HCV-spatially homogeneous infected equilibrium}
The objective of this section is to discuss the global
stability of the spatially homogeneous infected equilibrium $E^{*}$
for the PDE system (\ref{equation1.1}). We address global stability
by using the method of construction of Lyapunov functionals. These
Lyapunov functional is obtained from those for differential
equations by applying the method of Hattaf and Yousfi presented in
\cite{hattafyousfi2013}. We address this study with certain
assumptions namely : $u=0$ (i.e., there is no absorption effect),
$\alpha_{0}=1$ et $\alpha_{3}=\alpha_{1}\alpha_{2}$. Thus we have
the following results.
\begin{theorem}\label{theo3.5.1}
If $\mathcal{R}_{0} > 1$ the spatially homogeneous infected
equilibrium $E^{*}$ of PDE-model system (\ref{equation1.1}) is
globally asymptotically stable.
\end{theorem}
\begin{prooof}We first define the function
 \begin{eqnarray*}
G_{2}(H,I,V)&=& H-H^{*}-\int_{H^{*}}^{H}\frac{(\alpha + \rho)I^{*}}{\frac{(1-\eta)\beta \tau V^{*}}{(1+\alpha_{1}\tau )(1+\alpha_{2}V^{*})}} d\tau + I-I^{*}-I^{*}\ln(\frac{I}{I^{*}})\\
        &+& \frac{\alpha + \rho}{(1-\varepsilon)k}(\alpha_{0}+\alpha_{2}V^{*})\left( V-V^{*}-\int_{V^{*}}^{V}\frac{(\alpha + \rho)I^{*}}{\frac{(1-\eta)\beta\tau H^{*}}{(1+\alpha_{1}H^{*} )(1+\alpha_{2}\tau)}} d\tau \right).
\end{eqnarray*}
Then, the computation of the  derivative of $G_{2}$ with respect to
$t$ yields :

\begin{align*}
\frac{dG_{2}}{dt}=& \left[\lambda - dH -\alpha I -\frac{(\alpha+\rho)\mu}{(1-\varepsilon)k}V\right]\\
   &-(\alpha+\rho)I^{*}\frac{(1+\alpha_{1}H )(1+\alpha_{2} V^{*})}{(1-\eta)\beta HV^{*}}\left[ \lambda-dH-\frac{(1-\eta)\beta HV}{(1+\alpha_{1}H)(1+\alpha_{2}V)}+\rho I \right]\\
   & -\frac{I^{*}}{I}\left[\frac{(1-\eta)\beta HV}{(1+\alpha_{1}H)(1+\alpha_{2}V)}-(\alpha+\rho)I \right] -\frac{\alpha+\rho}{(1-\varepsilon)k}\frac{V^{*}}{V}\left[(1-\varepsilon)kI-\mu V \right].
\end{align*}
since
 $$\frac{(1-\eta)\beta
H^{*}V^{*}}{(1+\alpha_{1}H^{*})(1+\alpha_{2}V^{*})}=(\alpha+\rho)I^{*}$$
$$\lambda=dH^{*}+\alpha I^{*}$$
and,
$$\frac{(\alpha+\rho)\mu}{(1-\varepsilon)k}=\frac{(\alpha+\rho)I^{*}}{V^{*}}.$$
We have
\begin{eqnarray*}
 \frac{dG_{2}}{dt}  &=& \left[dH^{*}+\alpha I^{*} - dH -\alpha I -(\alpha+\rho)I^{*}
 \frac{V}{V^{*}}\right] \\
   &&-(\alpha+\rho)I^{*}\frac{(1+\alpha_{1}H )(1+\alpha_{2} V^{*})}
                {(1-\eta)\beta HV^{*}}\left[ dH^{*}+\alpha I^{*}-dH-\frac{(1-\eta)
                \beta HV}{(1+\alpha_{1}H)(1+\alpha_{2}V)}+\rho I \right]  \\
   &&-\frac{I^{*}}{I}\left[\frac{(1-\eta)\frac{ (1+\alpha_{1}H^{*})
                (1+\alpha_{2}V^{*})(\alpha+\rho)I^{*} }{(1-\eta) H^{*}V^{*}} HV}
                {(1+\alpha_{1}H)(1+\alpha_{2}V)}-(\alpha+\rho)I \right]
                -(\alpha+\rho)I\frac{V^{*}}{V} +\frac{(\alpha+\rho)\mu}{(1-\varepsilon)k}V^{*},
\end{eqnarray*}
\begin{align*}
                 =&\left[dH^{*}+(\alpha+\rho) I^{*}-\rho I^{*} - dH -\alpha I -(\alpha+\rho)I^{*}\frac{V}{V^{*}}\right] \\
                &-\Big[ \frac{H^{*}}{H} \frac{ 1+\alpha_{1}H}{1+\alpha_{1}H^{*}}dH^{*}+ \frac{H^{*}}{H} \frac{ 1+\alpha_{1}H}{1+\alpha_{1}H^{*}}\alpha I^{*}- \frac{ 1+\alpha_{1}H}{1+\alpha_{1}H^{*}}dH^{*}-\frac{V}{V^{*}}\frac{1+\alpha_{2}V^{*}}{1+\alpha_{2}V}(\alpha+\rho)I^{*}\\
                &+ \frac{H^{*}}{H}\frac{1+\alpha_{1}H}{1+\alpha_{1}H^{*}}\rho I \Big]\\
                &+ (\alpha+\rho)I^{*}\left[1-\frac{HI^{*}V(1+\alpha_{1}H^{*})
                (1+\alpha_{2}V^{*})}{H^{*}IV^{*}(1+\alpha_{1}H)(1+\alpha_{2}V)}
                \right] +(\alpha+\rho)I^{*}\left(1- \frac{I}{I^{*}}\frac{V^{*}}{V}\right),\\
\end{align*}
\begin{align*}
                =&dH^{*}\left[1-\frac{H}{H^{*}}-\frac{H^{*}}{H}\frac{1+\alpha_{1}H}{1+\alpha_{1}H^{*}}
                +\frac{1+\alpha_{1}H}{1+\alpha_{1}H^{*}}\right]\\
                &+ (\alpha+\rho)I^{*}\left[1-\frac{HI^{*}V(1+\alpha_{1}H^{*})
                (1+\alpha_{2}V^{*})}{H^{*}IV^{*}(1+\alpha_{1}H)(1+\alpha_{2}V)}
                +\frac{V}{V^{*}}\frac{1+\alpha_{2}V^{*}}{1+\alpha_{2}V} \right]\\
                &+(\alpha+\rho)I^{*}\left(2- \frac{V}{V^{*}} -\frac{I}{I^{*}}\frac{V^{*}}{V}\right)
                - \alpha I^{*}\left(\frac{\rho}{\alpha}+\frac{I}{I^{*}}\right)
                - \frac{H^{*}}{H}\frac{1+\alpha_{1}H}{1+\alpha_{1}H^{*}}\alpha I^{*}
                -\frac{H^{*}}{H}\frac{1+\alpha_{1}H}{1+\alpha_{1}H^{*}}\rho I,\\
\end{align*}
\begin{align*}
                =&-\frac{d(H-H^{*})^{2}}{H(1+\alpha_{1}H^{*})} + (\alpha+\rho)I^{*}\left[-1-\frac{V}{V^{*}}  +\frac{V}{V^{*}}\frac{1+\alpha_{2}V^{*}}{1+\alpha_{2}V} + \frac{1+\alpha_{2}V}{1+\alpha_{2}V^{*}}\right]\\
                &+ (\alpha+\rho)I^{*}\left[4- \frac{H^{*}}{H}\frac{1+\alpha_{1}H}{1+\alpha_{1}H^{*}}- \frac{HI^{*}V(1+\alpha_{1}H^{*})(1+\alpha_{2}V^{*})}{ H^{*}IV^{*}(1+\alpha_{1}H)(1+\alpha_{2}V)} -\frac{I}{I^{*}}\frac{V^{*}}{V}-\frac{1+\alpha_{2}V}{1+\alpha_{2}V^{*}}\right]\\
                & - \alpha I^{*}\left(\frac{\rho}{\alpha}+\frac{I}{I^{*}}\right) - \frac{H^{*}}{H}\frac{1+\alpha_{1}H}{1+\alpha_{1}H^{*}}\alpha I^{*} -\frac{H^{*}}{H}\frac{1+\alpha_{1}H}{1+\alpha_{1}H^{*}}\rho I +
                (\alpha+\rho)I^{*}\frac{H^{*}}{H}\frac{1+\alpha_{1}H}{1+\alpha_{1}H^{*}}.
\end{align*}
Therefore
\begin{align*}
\frac{dG_{2}}{dt}=&-\frac{d(H-H^{*})^{2}}{H(1+\alpha_{1}H^{*})}-\frac{\alpha_{2}(\alpha+\rho)
I^{*}(V-V^{*})^{2}}{V^{*}(1+\alpha_{2}V^{*})(1+\alpha_{2}V)}-
\alpha I^{*}\left(\frac{\rho}{\alpha}+\frac{I}{I^{*}}\right) \\
                &+ (\alpha+\rho)I^{*}\left[4- \frac{H^{*}}{H}\frac{1+\alpha_{1}H}
                {1+\alpha_{1}H^{*}}- \frac{HI^{*}V(1+\alpha_{1}H^{*})(1+\alpha_{2}V^{*})}
                { H^{*}IV^{*}(1+\alpha_{1}H)(1+\alpha_{2}V)} -\frac{I}{I^{*}}\frac{V^{*}}{V}
                -\frac{1+\alpha_{2}V}{1+\alpha_{2}V^{*}}\right].
\end{align*}
 We have :
  $$ (\alpha+\rho)I^{*}\left[4-
\frac{H^{*}}{H}\frac{1+\alpha_{1}H}{1+\alpha_{1}H^{*}}-
\frac{HI^{*}V(1+\alpha_{1}H^{*})(1+\alpha_{2}V^{*})}{
H^{*}IV^{*}(1+\alpha_{1}H)(1+\alpha_{2}V)}
-\frac{I}{I^{*}}\frac{V^{*}}{V}-\frac{1+\alpha_{2}V}{1+\alpha_{2}V^{*}}\right]\leq
0$$
 since the left side of the latter inequality is the difference
between the geometric mean and the arithmetic mean. That is
$\frac{dG_{2}}{dt}\leq 0$. Otherwise $\frac{dG_{2}}{dt}=0$ if and
only if $H=H^{*}$, $I=I^{*}$ et $V=V^{*}$. Thus $G_{2}$ is a
Lyapunov functional of the differential equation associated to the
PDE-model system (\ref{equation1.1}). Therefore using
Lyapunov-LaSalle invariance principle \cite{khalil1} combined to the
method presented in \cite{hattafyousfi2013}, the  functional defined
by
$$L_{2}=\int_{\Omega} G_{2}dx $$
is a lyapunov functional of the PDE-model system (\ref{equation1.1})
at the spatially homogeneous infected equilibrium $E^{*}$. Therefore
$E^{*}$ is  globally  asymptotically stable. This completes the
proof of Theorem~\ref{theo3.5.1}.
\end{prooof}
\section{Numerical simulations}\label{sec:4}
\begin{figure}[!h]
\centering
\begin{subfigure}[]{}
\includegraphics[angle=0,height=6cm,width
=5.5cm]{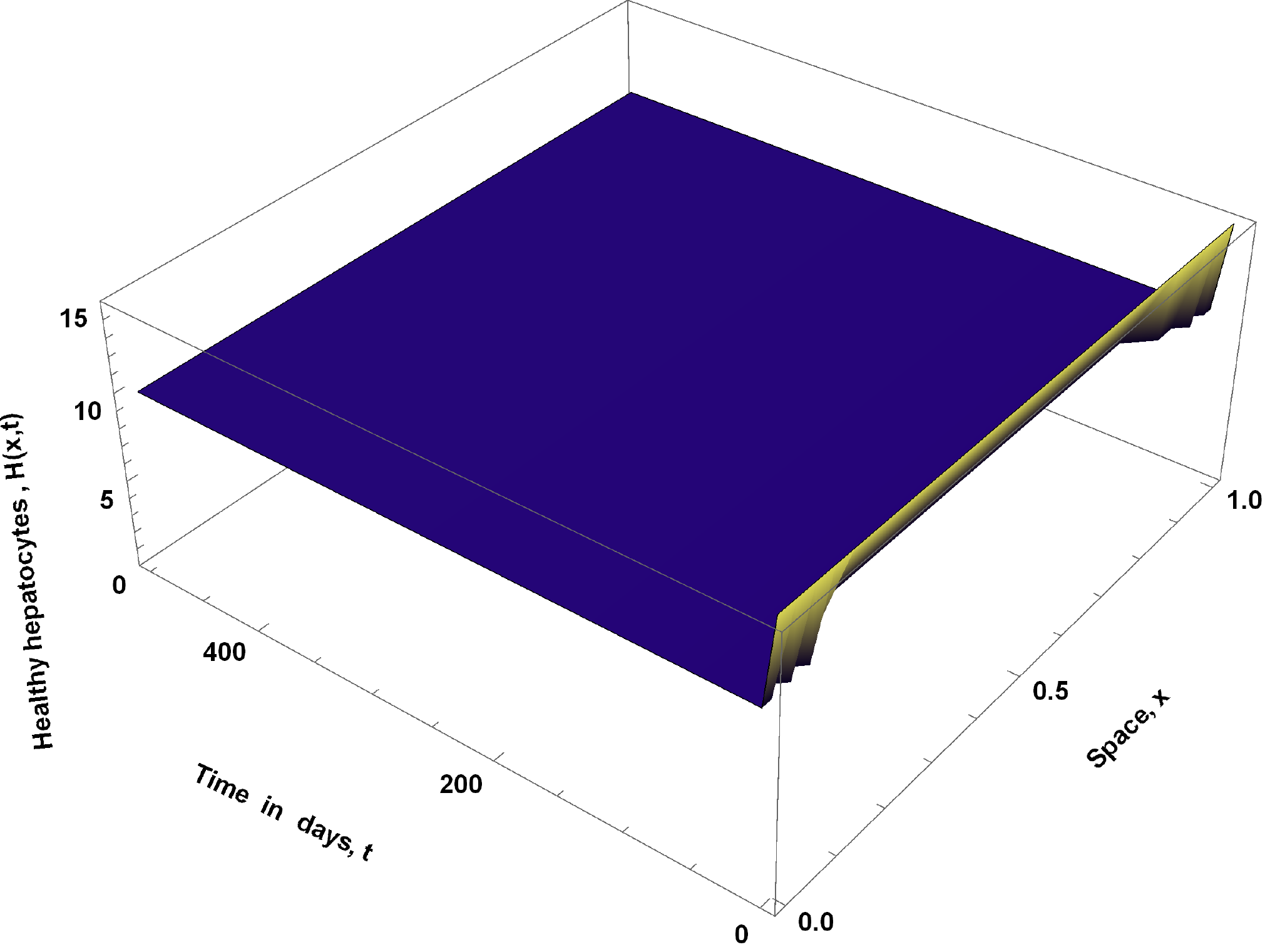}
\end{subfigure}
\begin{subfigure}[]{}
\includegraphics[angle=0,height=6cm,
width=5.5cm]{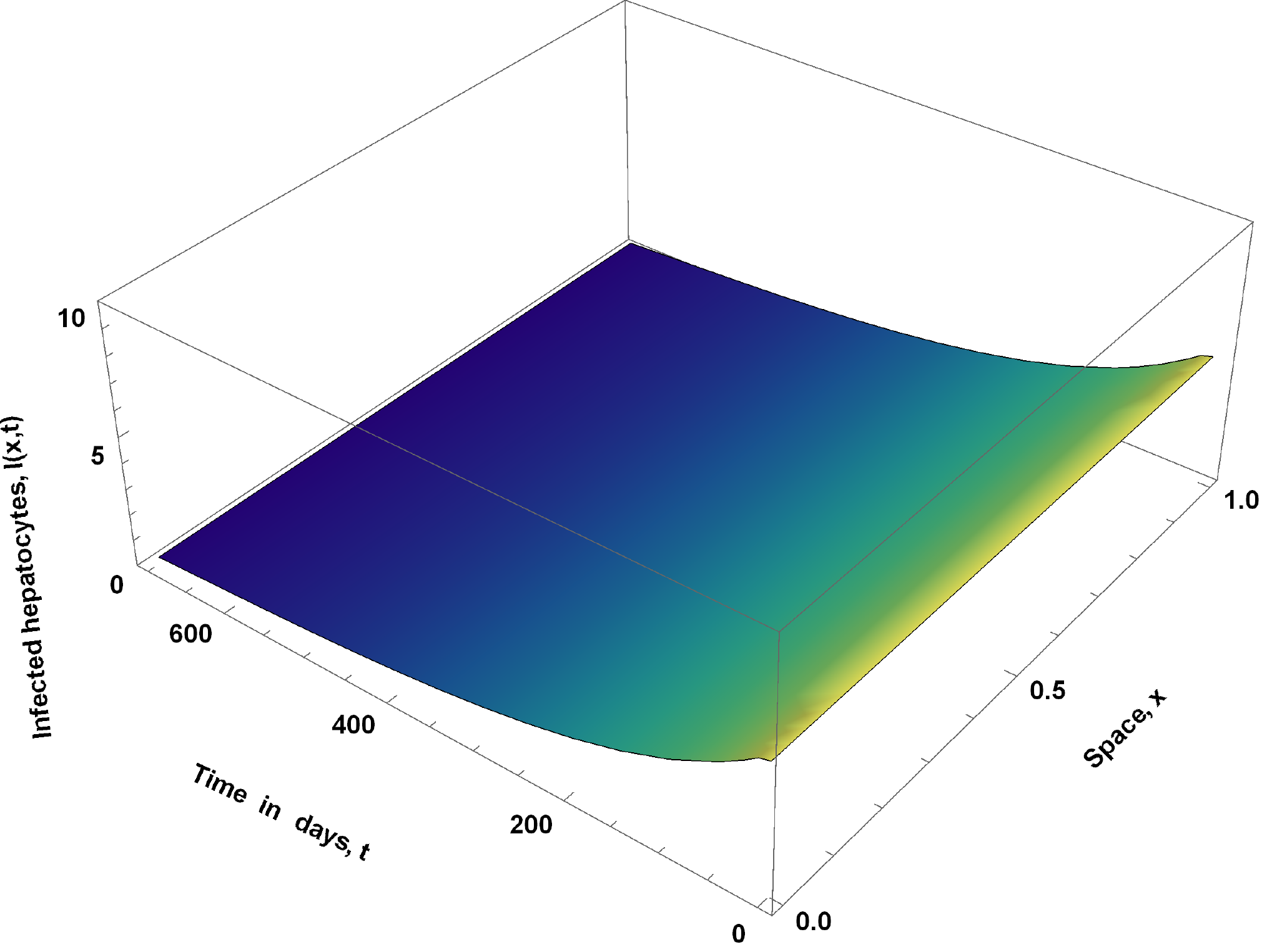}
\end{subfigure}
\begin{subfigure}[]{}
\includegraphics[angle=0,height=6cm,width=7cm]
{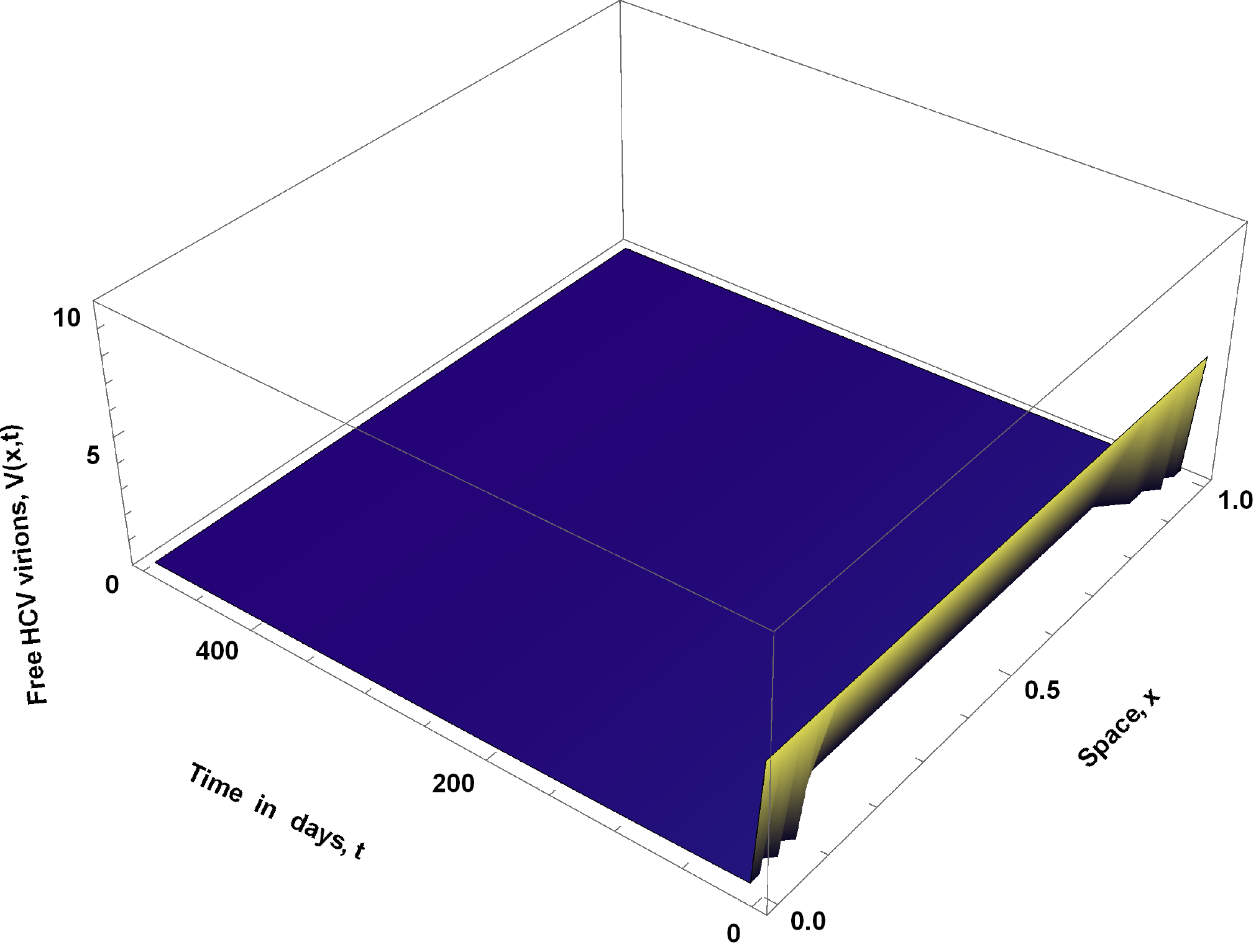}
\end{subfigure}
\caption{Simulations of IBVP (\ref{equation1.1}) under Neumann
boundary conditions (\ref{1.5}) and initial condition (\ref{4.5a})}
\label{figg1}
\end{figure}

In this section, we present the numerical simulations to illustrate
our theoretical results. To simplify, we consider IBVP
(\ref{equation1.1}) with $\Omega = (1)$ under Neumann boundary
condition
\begin{equation}\label{1.5}
    \frac{\partial H}{\partial \nu} =0, \;\frac{\partial I}{\partial \nu} =0, \;\frac{\partial V}{\partial \nu} =0 \;\; t>0, \;\;x= 1
\end{equation}
and, following initial conditions
\begin{equation}\label{4.5a}
  H(x, 0)=5,       \;\;   I(x, 0)=5,     \;\; V(x, 0)=5,
\end{equation}
and
\begin{equation}\label{4.5b}
  H(x, 0)=15,      \;\;   I(x, 0)=5,     \;\; V(x, 0)=5.
\end{equation}
Now we choose the  numerical values of the parameters for the
PDE-cellular model system (\ref{equation1.1}) as follows:
$\lambda=50$; $d=5$; $\rho=0,01$; $\alpha=0,05$; $D_{1}=D_{2}=D_{3}=0,1$;
$\eta=0,00004$;
 $\alpha_{3}= 0,03$; $\varepsilon=0,5$; $\alpha_{2}= 0,02$;
 $k=2$; $\alpha_{1}=0,1$; $\mu=20$; $\alpha_{0}=1$; $\beta=0,24$ et
 $u=1$. By calculation we have $\mathcal{R}_{0}=0.943361$. In this
 case, PDE-cellular model system (\ref{equation1.1}) has a spatially
 homogeneous equilibrium $E_{0}=(10, 0, 0)$. Hence by Theorem \ref{theogloE0}
$E_{0}$ is globally asymptotically stable. Numerical simulation
illustrates our result (see figure \ref{figg1}).

\begin{figure}[!h]
\centering
\begin{subfigure}[]{}
\includegraphics[angle=0,height=6.2cm,width
=5.5cm]{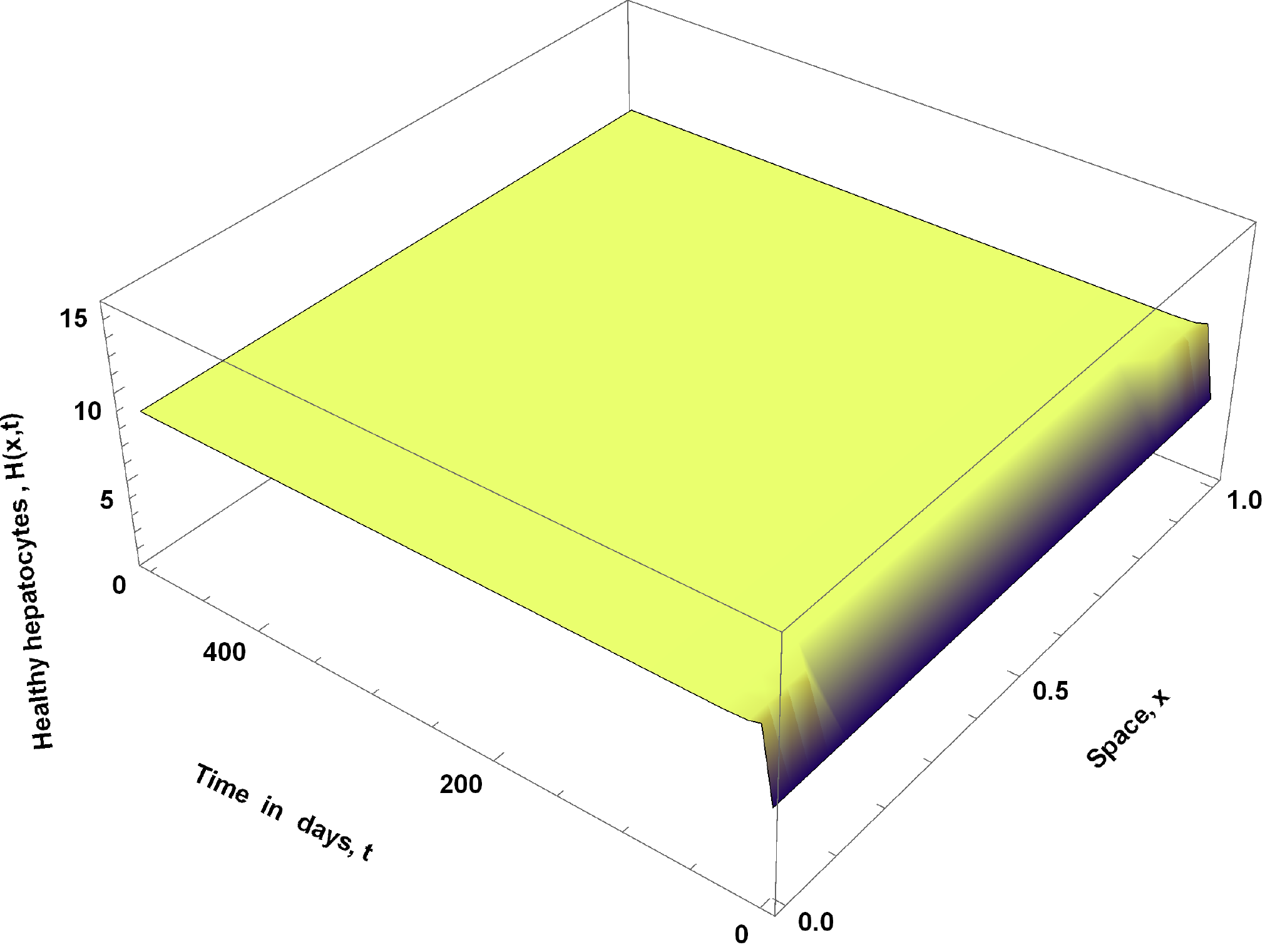}
\end{subfigure}
\begin{subfigure}[]{}
\includegraphics[angle=0,height=6.2cm,
width=5.5cm]{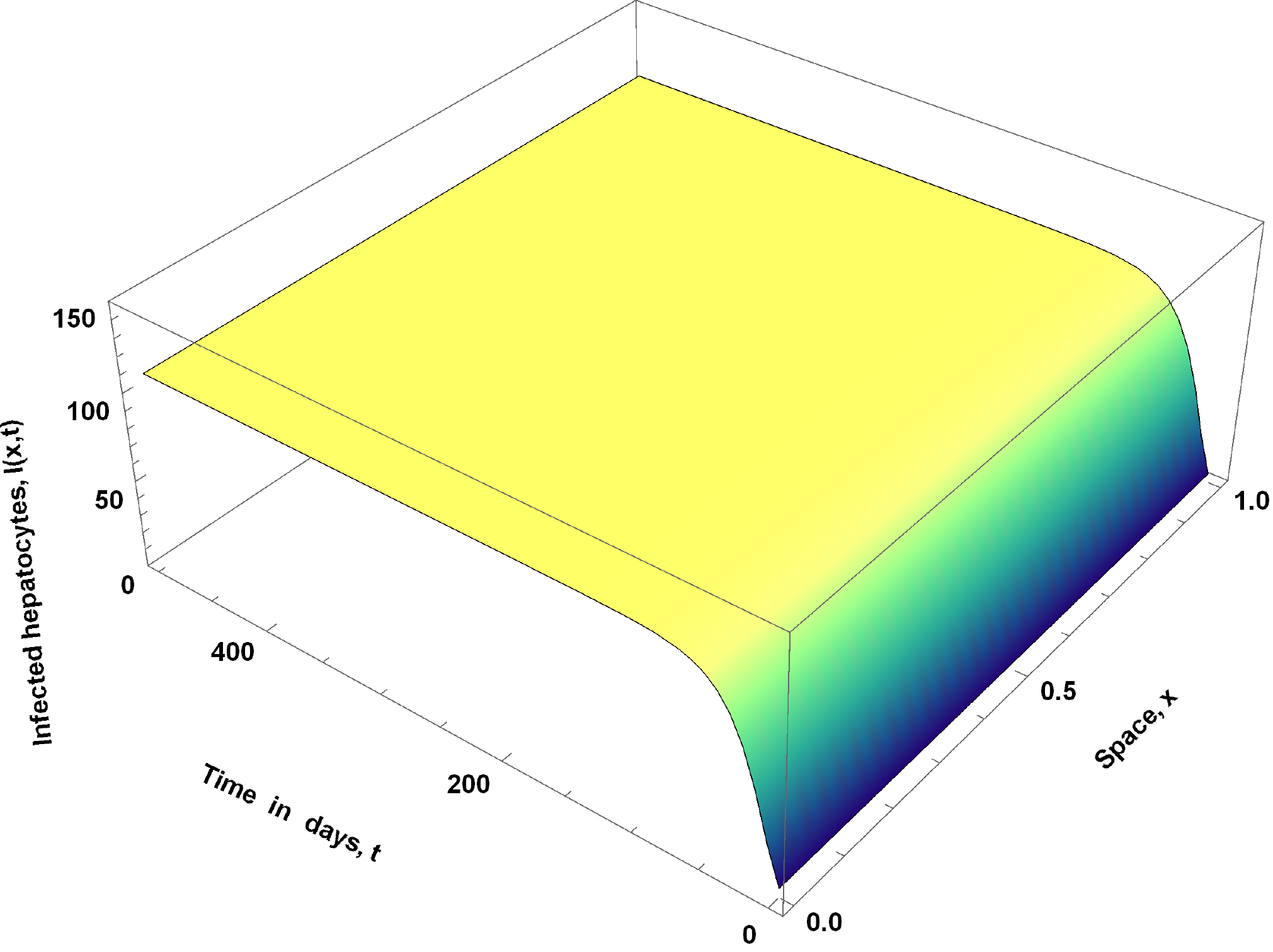}
\end{subfigure}
\begin{subfigure}[]{}
\includegraphics[angle=0,height=6cm,width=7cm]
{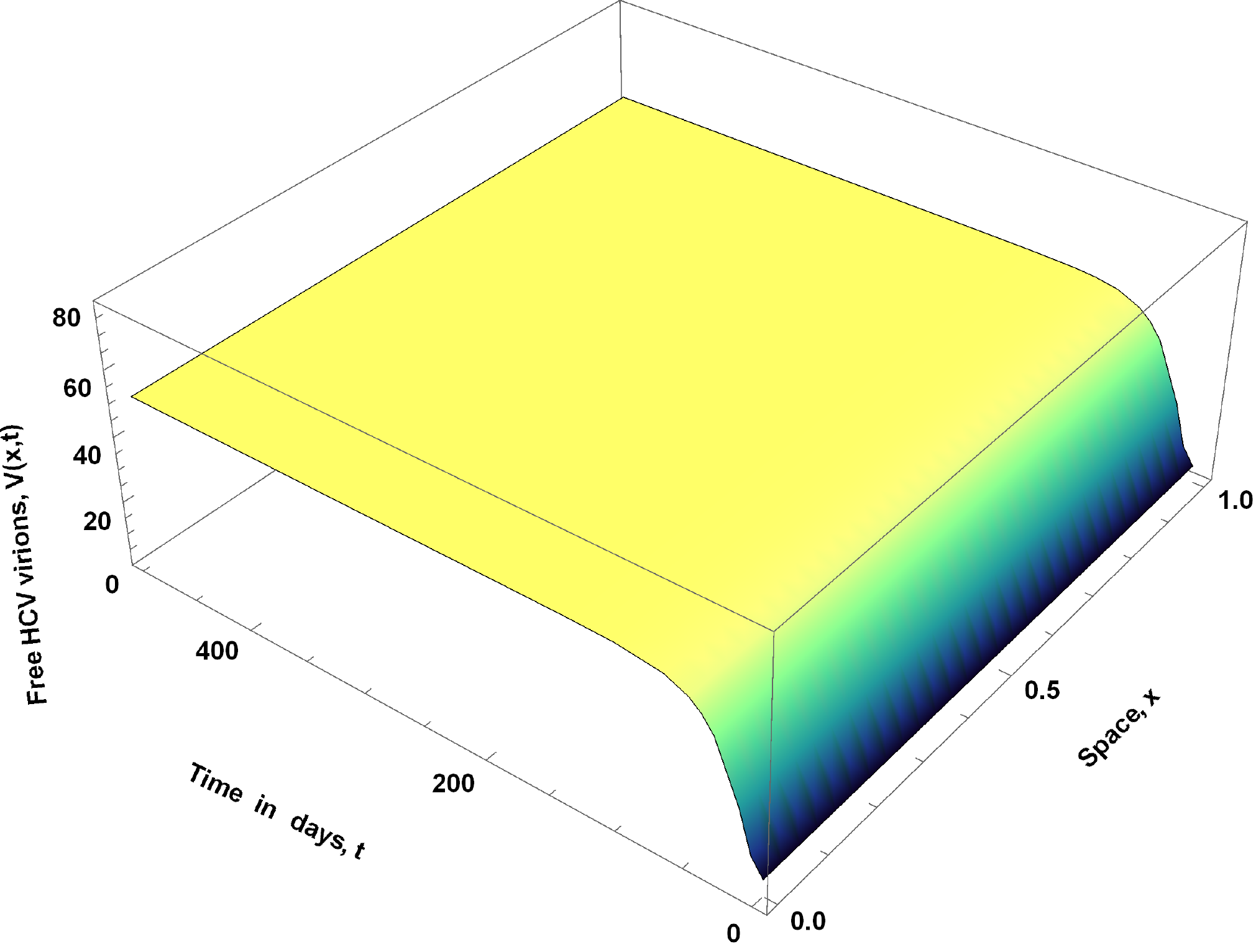}
\end{subfigure}
\caption{Simulations of IBVP (\ref{equation1.1}) under Neumann
boundary conditions (\ref{1.5}) and initial condition (\ref{4.5b})}
\label{figg2}
\end{figure}

Otherwise we choose the  numerical values of the parameters for the
PDE-cellular model system (\ref{equation1.1}) as follows :
$\lambda=50$; $d=5$; $\rho=0,01$; $\alpha=0,05$; $D_{1}=D_{2}=D_{3}=0,1$;
$\eta=0,00004$; $\varepsilon=0,5$;
  $\alpha_{0}=1$; $\alpha_{1}=0,1$; $\alpha_{2}= 0,02$; $\alpha_{3}= 0,03$; $k=2$; $\mu=2$;
  $\beta=0,24$ et $u=1$. By calculation we have $\mathcal{R}_{0}=6.25009$. In this
 case, PDE-cellular model system (\ref{equation1.1}) has a spatially
 homogeneous equilibrium $E^{*}=(5; 500; 235)$. Hence by Theorem \ref{theo3.5.1}
$E^{*}$ is globally asymptotically stable. Numerical simulation
illustrates our result (see figure \ref{figg2}).
\begin{acknowledgements}
The authors would like to thank the editor and the anonymous
referees for their valuable remarks and comments which have led to
the improvement of the quality of their study.
\end{acknowledgements}
\section*{Conflict of interest}
 The authors declare that they have no conflict
 of interests regarding the publication of this paper.
\section{Conclusion}\label{sec:5}
In this work, we addressed the dynamics of a reaction diffusion HCV
intra-host infection model with the Hattaf-Yousfi incidence rate,
which is a generalized nonlinear incidence rate. The object of this
work was to make a mathematical analysis of a cellular model of HCV
infection which assumes that virions diffuse into the liver, which
uses the Hattaf-Yousfi functional response generalizing most of the
functional responses that exist. Our model also takes into account
the absorption effect which is much neglected in the literature. We
first showed that the initial value and boundary problem
(\ref{equation1.1}) admits a unique global solution in time.
 And secondly, we have shown that this unique solution is positive and
 uniformly bounded. Then, we gave the expression of
 basic reproduction number $\mathcal{R}_{0}$ which
 is the parameter from which we studied the dynamics of our model at
 equilibria whose existence and
 uniqueness of these have been previously proven. More precisely,
 we have shown that, if $ \mathcal{R}_{0} <1 $, the unique uninfected
 equilibrium
  point is locally and globally asymptotically stable. This means that,
  under this condition,infection disappears.
Otherwise, the uninfected equilibrium point is unstable; and in this
case the infection persists in the host. It has also been shown
under the hypothesis $ \mathcal{R}_{0}> 1 $, that the infected
equilibrium point is locally and globally asymptotically stable.
This edifying work ended with numerical simulations, carried out on
the Mathematica software, which confirmed our theoretical results.
 In order to get as close as possible to complex reality
of biological phenomena, we envisage in the future, the mathematical
analysis of models taking into account cell proliferation and the delay.
\appendix
\section{Proof of Proposition~\ref{propo3.9}}
The proof is established by using Banach's Fixed Point Theorem.\\
  Choose $\beta$ such that $\frac{3}{4}<\beta<1$, then the injection
   $ \mathcal{I}: D(\mathcal{H}^{\beta})\rightarrow C_{B}^{0}$
   is continuous by Lemma~\ref{lem2.1.2}.
   For $r_{0}>0$ and $T>0$, the following closed ball is considered
   :
  \begin{equation}\label{ball1}
  B_{r_{0}}(H_{0},I_{0},V_{0})=\left\{(H,I,V)\in(C_{B}^{0}(]0,T],
  D(\mathcal{H}^{\beta})))^{3}: \|H-H_{0}\|_{\beta},
  \|I-I_{0}\|_{\beta},\|V-V_{0}\|_{\beta}\leq r_{0} \right\}.
  \end{equation}
  By Proposition~\ref{proposition2.1.1},
  we have the local Lipschitz properties:
\begin{eqnarray}
\nonumber
\|\mathcal{F}(t,H_{1},I_{1},V_{1})-\mathcal{F}(t,H_{2},I_{2},V_{2})\|_{2}
&\leq&
K_{1}^{1}\|H_{1}-H_{2}\|_{D(\mathcal{H}^{\beta})}+K_{2}^{1}\|I_{1}
-I_{2}\|_{D(\mathcal{H}^{\beta})}+K_{3}^{1}\|V_{1}-V_{2}\|_{D(\mathcal{H}^{\beta})},  \\
\nonumber\|\mathcal{G}(t,H_{1},I_{1},V_{1})-\mathcal{G}(t,H_{2},I_{2},V_{2})\|_{2}
   &\leq& K_{1}^{2}\|H_{1}-H_{2}\|_{D(\mathcal{H}^{\beta})}+K_{2}^{2}
   \|I_{1}-I_{2}\|_{D(\mathcal{H}^{\beta})}
   +K_{3}^{2}\|V_{1}-V_{2}\|_{D(\mathcal{H}^{\beta})},
    \\
\nonumber\|\mathcal{Q}(t,H_{1},I_{1},V_{1})-\mathcal{Q}(t,H_{2},I_{2},V_{2})\|_{2}
&\leq&
K_{1}^{3}\|H_{1}-H_{2}\|_{D(\mathcal{H}^{\beta})}+K_{2}^{3}\|I_{1}
-I_{2}\|_{D(\mathcal{H}^{\beta})}+K_{3}^{3}\|V_{1}
-V_{2}\|_{D(\mathcal{H}^{\beta})},\\
&& \label{3.20}
\end{eqnarray}
for $t\in [0,T]$, $(H_{1},I_{1},V_{1}), (H_{2},I_{2},V_{2}) \in
B_{r_{0}}(H_{0},I_{0},V_{0})$ and Lipschitz-constants $ K_{j}^{i}
>0$, $i,j=1,2,3$. In addition $(y_{1},y_{2},y_{3})\in
B_{r_{0}}(H_{0},I_{0},V_{0})$, define $P:[0,T]\rightarrow
L^{2}(\Omega)$, $Q:[0,T]\rightarrow L^{2}(\Omega)$ and
$R:[0,T]\rightarrow L^{2}(\Omega)$ as follows
 \begin{eqnarray*}
  Py_{1}(t) &=& G_{1}(t)H_{0}+ \int_{0}^{t} G_{1}(t-\tau)\mathcal{F}(t,y_{1}(\tau),y_{2}(\tau),y_{3}(\tau))d\tau, \\
  Qy_{2}(t) &=& G_{2}(t)I_{0} + \int_{0}^{t} G_{2}(t-\tau)\mathcal{G}(t,y_{1}(\tau),y_{2}(\tau),y_{3}(\tau))d\tau, \\
  Ry_{3}(t) &=& G_{3}(t)V_{0} + \int_{0}^{t} G_{3}(t-\tau)
  \mathcal{Q}(t,y_{1}(\tau),y_{2}(\tau),y_{3}(\tau))d\tau.
  \end{eqnarray*}
 Finally, set
 $M_{1}=\sup\limits_{t\in[0,T]} \|\mathcal{F}
 (t,y_{1}(0),y_{2}(0),y_{3}(0))\|_{2}$,
 $M_{2}=\sup\limits_{t\in[0,T]}
  \|\mathcal{G}(t,y_{1}(0),y_{2}(0),y_{3}(0))\|_{2}$,
  $M_{3}=\sup\limits_{t\in[0,T]}\|\mathcal{Q}(t,y_{1}(0),y_{2}(0),y_{3}(0))\|_{2}$
 and choose T so that:
\begin{equation}\label{equation2.10}
\|G_{1}(h)H_{0}-H_{0}\|_{\beta},\|G_{2}(h)I_{0}-I_{0}\|_{\beta},
 \|G_{3}(h)V_{0}-V_{0}\|_{\beta}\leq \frac{3r_{0}}{4},
 \;\; 0\leq h\leq T,
\end{equation}
\begin{equation}\label{equation2.11}
C_{\beta,2}^{i} (M_{i} + r_{0}K_{1}^{i} + r_{0}K_{2}^{i}
  + r_{0}K_{3}^{i})\int_{0}^{T} s^{-\beta}ds \leq\frac{r_{0}}{4},
   \;\; i=1,2,3,
\end{equation}
where $C_{\beta,2}^{i}$ is the constant in property 2)
  of Corollary~\ref{corol2.1.1} for operator A ($A$ is zero for
  $i=1$, $2$). Note that such T
exists since $G_{1}(h)$, $G_{2}(h)$ and $G_{3}(h)$ converge to $Id$
as $h$ tends to $0^{+}$ by definition of an analytic semigroup and
  $$\int_{0}^{T} s^{-\beta}ds=\frac{1}{1-\beta}h^{1-\beta}
  \rightarrow 0, \;\; h\rightarrow 0^{+} \;\;for \;\; \beta<1.$$
  Then the proof will continue according to the following two points
  :
  \begin{description}
    \item[(a)] it is shown that (P,Q,R) maps
    $B_{r_{0}}(H_{0},I_{0},V_{0})$ into itself,
    \item[(b)] it is shown that (P,Q,R) is a strict
    contraction on $B_{r_{0}}(H_{0},I_{0},V_{0})$,
    allowing the use of Banach's Fixed
Point Theorem to get the existence of a unique fixed point in
$B_{r_{0}}(H_{0},I_{0},V_{0})$.
  \end{description}
  Let $(H,I,V)\in B_{r_{0}}(H_{0},I_{0},V_{0})$.
  Then, using (\ref{equation2.10}) and property 2)
  of Corollary \ref{corol2.1.1}, we have :
  \begin{eqnarray*}
     \|PH(t)-H_{0} \|_{D(\mathcal{H}^{\beta})}&=& \Big\|G_{1}(t)H_{0}-H_{0}+ \int_{0}^{t} G_{1}(t-\tau)\mathcal{F}(\tau,H(\tau),
     I(\tau),V(\tau))d\tau \Big\|_{D(\mathcal{H}^{\beta})} \nonumber, \\
      & \leq & \|G_{1}(t)H_{0}-H_{0}\|_{D(\mathcal{H}^{\beta})}
      + \int_{0}^{t} \|G_{1}(t-\tau)\mathcal{F}(\tau,H(\tau),I(\tau),V(\tau))
      \|_{D(\mathcal{H}^{\beta})}d\tau, \nonumber\\
      & \leq & \frac{3r_{0}}{4} + \int_{0}^{t} C_{\beta,2}^{1}(t-s)^{-\beta}
       \Big\|\mathcal{F}(\tau,H(\tau),I(\tau),V(\tau)) \Big\|_{2}d\tau, \nonumber \\
       & \leq & \frac{3r_{0}}{4} + \int_{0}^{t} C_{\beta,2}^{1}(t-s)^{-\beta}
       \Big\|\mathcal{F}(\tau,H(\tau),I(\tau),V(\tau)) -\mathcal{F}(\tau,H_{0},I_{0},V_{0})
       + \mathcal{F}(\tau,H_{0},I_{0},V_{0})\Big\|_{2}d\tau, \nonumber \\
       & \leq & \frac{3r_{0}}{4} + \int_{0}^{t} C_{\beta,2}^{1}(t-s)^{-\beta} \Big
       (\Big\|\mathcal{F}(\tau,H(\tau),I(\tau),V(\tau)) -\mathcal{F}(\tau,H_{0},I_{0},V_{0})\Big\|_{2}\\
       &+& \Big\|\mathcal{F}(\tau,H_{0},I_{0},V_{0})\Big\|_{2}\Big)d\tau, \nonumber
\end{eqnarray*}
\begin{eqnarray*}
      \|PH(t)-H_{0} \|_{D(\mathcal{H}^{\beta})}& \leq & \frac{3r_{0}}{4} + \int_{0}^{t} C_{\beta,2}^{1}(t-s)^{-\beta}
       \Big(K_{1}^{1}r_{0} + K_{2}^{1}r_{0} + K_{3}^{1}r_{0} + M_{1}\Big)d\tau, \nonumber\\
      & \leq & \frac{3r_{0}}{4} + C_{\beta,2}^{1}\left(K_{1}^{1}r_{0} + K_{2}^{1}r_{0}
        + K_{3}^{1}r_{0} + M_{1}\right) \int_{0}^{t} (t-s)^{-\beta}d\tau, \nonumber\\
      & \leq & \frac{3r_{0}}{4}+ 0,\\
      & \leq & r_{0} \;\; \mbox{for} \;\; 0\leq t\leq T,
  \end{eqnarray*}
 and similarly
 $$\|QI(t)-I_{0} \|_{D(\mathcal{H}^{\beta})} \leq  r_{0} \;\;
  \mbox{and} \;\;<
 \|RV(t)-V_{0} \|_{D(\mathcal{H}^{\beta})} \leq  r_{0}.$$
  Showing that $(P,Q,R)$ maps $B_{r_{0}}(H_{0},I_{0},V_{0})$
into itself. Furthermore, from property 1) in
Corollary~\ref{corol2.1.1} and Lemma~\ref{lem2.1.3} we compute:
$$\|PH(t+h)-PH(t) \|_{D(\mathcal{H}^{\beta})}  $$
\begin{eqnarray*}
 &=& \Big\|G_{1}(t+h)H_{0}-G_{1}(t)H_{0}
 + \int_{0}^{t+h} G_{1}(t+h-\tau)\mathcal{F}(\tau,H(\tau),I(\tau),V(\tau))d\tau\\
 &&- \int_{0}^{t} G_{1}(t-\tau)\mathcal{F}(\tau,H(\tau),I(\tau),V(\tau))d\tau
 \Big\|_{D(\mathcal{H}^{\beta})},\\
 &=& \Big\|G_{1}(t)G_{1}(h)H_{0}-G_{1}(t)H_{0}+ \int_{0}^{t} G_{1}(h)G_{1}(t-\tau)
 \mathcal{F}(\tau,H(\tau),I(\tau),V(\tau))d\tau \\
  & &+ \int_{t}^{t+h} G_{1}(h)G_{1}(t-\tau)\mathcal{F}(\tau,H(\tau),I(\tau),V(\tau))
  d\tau - \int_{0}^{t} G_{1}(t-\tau)\mathcal{F}(\tau,H(\tau),I(\tau),V(\tau))d\tau
  \Big\|_{D(\mathcal{H}^{\beta})},\\
  &=& \Big\|G_{1}(t)(G_{1}(h)-Id)H_{0}+ \int_{0}^{t} (G_{1}(h)-Id)G_{1}
  (t-\tau)\mathcal{F}(\tau,H(\tau),I(\tau),V(\tau))d\tau \\
 & &+ \int_{t}^{t+h} G_{1}(t+h-\tau)\mathcal{F}(\tau,H(\tau),I(\tau),V(\tau))d\tau
  \Big\|_{D(\mathcal{H}^{\beta})}.
  \end{eqnarray*}
The strong continuity of the semigroup and Lemma~\ref{lem2.1.3}
yields
 $$\|PH(t+h)-PH(t) \|_{D(\mathcal{H}^{\beta})}\rightarrow 0 \;\;
\mbox{as} \;\; h\rightarrow 0^{+}.$$
Therefore P is continuous from [0,T] into $D(\mathcal{H}^{\beta})$.\\
A similar Calculation shows that Q and R have the same properties
and item (a) is proved.
\\\noindent
Presently let us show that $(P,Q,R)$ is a strict contraction on $B_{r_{0}}(H_{0},I_{0},V_{0})$.\\
Let $(y_{1},y_{2},y_{3}),\; (z_{1},z_{2},z_{3})\in
B_{r_{0}}(H_{0},I_{0},V_{0})$. Then
  \begin{eqnarray*}
 \|Py_{1}(t)-Pz_{1}(t) \|_{D(\mathcal{H}^{\beta})} &=&
  \Big\|\int_{0}^{t} G_{1}(t-\tau)\Big (\mathcal{F}(\tau,y_{1}
  (\tau),y_{2}(\tau),y_{3}(\tau))-\mathcal{F}(\tau,z_{1}(\tau),
  z_{2}(\tau),z_{3}(\tau))\Big)d\tau \Big\|_{D(\mathcal{H}^{\beta})},\\
 &\leq & \int_{0}^{t} \Big\|G_{1}(t-\tau)\Big (\mathcal{F}
 (\tau,y_{1}(\tau),y_{2}(\tau),y_{3}(\tau))-\mathcal{F}
 (\tau,z_{1}(\tau),z_{2}(\tau),z_{3}(\tau))\Big)
 \Big\|_{D(\mathcal{H}^{\beta})}d\tau,\\
  &\leq & \int_{0}^{t} C_{\beta,2}^{1}(t-\tau)^{-\beta}
  \Big\|\mathcal{F}(\tau,y_{1}(\tau),y_{2}(\tau),y_{3}(\tau))
  -\mathcal{F}(\tau,z_{1}(\tau),z_{2}(\tau),z_{3}(\tau))
  \Big\|_{2}d\tau,\\
   &\leq & \int_{0}^{t} C_{\beta,2}^{1}(t-\tau)^{-\beta}
   \Big( K_{1}^{1}\|y_{1}-z_{1}\|_{D(\mathcal{H}^{\beta})}
   +K_{2}^{1}\|y_{2}-z_{2}\|_{D(\mathcal{H}^{\beta})}
  +K_{3}^{1}\|y_{3}-z_{3}\|_{D(\mathcal{H}^{\beta})}\Big)
   d\tau,\\
   &\leq & C_{\beta,2}^{1} \Big( K_{1}^{1} \sup_{\tau\in [0,t]}
   \|y_{1}-z_{1}\|_{D(\mathcal{H}^{\beta})}+K_{2}^{1}
   \sup_{\tau\in [0,t]}\|y_{2}-z_{2}\|_{D(\mathcal{H}^{\beta})}\\
   &+& K_{3}^{1}\sup_{\tau\in [0,t]}
  \|y_{3}-z_{3}\|_{D(\mathcal{H}^{\beta})}\Big)
  \int_{0}^{t} (t-\tau)^{-\beta} d\tau.
  \end{eqnarray*}
That is :
\begin{equation*}
    \|Py_{1}(t)-Pz_{1}(t) \|_{D(\mathcal{H}^{\beta})}
\end{equation*}
\begin{eqnarray*}
&\leq & C_{\beta,2}^{1} \Big(
K_{1}^{1}+K_{2}^{1}+K_{3}^{1}+\frac{M_{1}} {r_{0}}\Big)
\Big(\sup_{\tau\in [0,t]} \|y_{1}-z_{1}
\|_{D(\mathcal{H}^{\beta})}+\sup_{\tau\in [0,t]}\|y_{2}-z_{2}
\|_{D(\mathcal{H}^{\beta})} \\
& &+\sup_{\tau\in [0,t]}\|y_{3}-z_{3}\|_{D(\mathcal{H}^{\beta})}
\Big) \int_{0}^{t} (t-\tau)^{-\beta} d\tau,\\
&\leq & \Big(\sup_{\tau\in [0,t]} \|y_{1}-z_{1}\|_{D(\mathcal{H}^{\beta})}+\sup_{\tau\in [0,t]}\|y_{2}-z_{2}\|_{D(\mathcal{H}^{\beta})}+\sup_{\tau\in [0,t]}\|y_{3}-z_{3}\|_{D(\mathcal{H}^{\beta})}\Big)\\
& &\times C_{\beta,2}^{1} \Big( K_{1}^{1}+K_{2}^{1}
+K_{3}^{1}+\frac{M_{1}}{r_{0}}\Big) \int_{0}^{t} s^{-\beta} ds,\\
&\leq & \Big(\sup_{\tau\in [0,t]} \|y_{1}-z_{1}\|_{D(\mathcal{H}^{\beta})}+\sup_{\tau\in [0,t]}\|y_{2}-z_{2}\|_{D(\mathcal{H}^{\beta})}+\sup_{\tau\in [0,t]}\|y_{3}-z_{3}\|_{D(\mathcal{H}^{\beta})}\Big)\\
& &\times \frac{1}{r_{0}}C_{\beta,2}^{1} \left(
r_{0}K_{1}^{1}+r_{0}K_{2}^{1}+r_{0}K_{3}^{1}+M_{1}\right)
\int_{0}^{T} s^{-\beta} ds,\\
&\leq & \Big(\sup_{\tau\in [0,t]} \|y_{1}-z_{1}\|_{D(\mathcal{H}^
{\beta})}+\sup_{\tau\in [0,t]}\|y_{2}-z_{2}\|
_{D(\mathcal{H}^{\beta})}+\sup_{\tau\in [0,t]}
\|y_{3}-z_{3}\|_{D(\mathcal{H}^{\beta})}\Big)
\frac{1}{r_{0}}\times\frac{r_{0}}{4},\\
&\leq & \frac{1}{4} \Big( \sup_{\tau\in [0,t]}
\|y_{1}-z_{1}\|_{D(\mathcal{H}^{\beta})} +\sup_{\tau\in
[0,t]}\|y_{2}-z_{2}\|_{D(\mathcal{H}^{\beta})} + \sup_{\tau\in
[0,t]}\|y_{3}-z_{3}\|_{D(\mathcal{H}^{\beta})}\Big),
\end{eqnarray*}

 for every $t\in [0,T]$. Hence
\begin{small}
\begin{eqnarray*}
\sup_{t\in [0,T]}\|Py_{1}(t)-Pz_{1}(t) \|_{D(\mathcal{H}^{\beta})}
&\leq& \frac{1}{4} \Big( \sup_{\tau\in [0,t]}
\|y_{1}-z_{1}\|_{D(\mathcal{H}^{\beta})} +\sup_{\tau\in
[0,t]}\|y_{2}-z_{2}\|_{D(\mathcal{H}^{\beta})} + \sup_{\tau\in
[0,t]}\|y_{3}-z_{3}\|_{D(\mathcal{H}^{\beta})}\Big).
\end{eqnarray*}
\end{small}
Similarly
\begin{small}
\begin{eqnarray*}
\sup_{t\in [0,T]}\|Qy_{1}(t)-Qz_{1}(t) \|_{D(\mathcal{H}^{\beta})}
&\leq &  \frac{1}{4} \Big( \sup_{\tau\in [0,t]}
\|y_{1}-z_{1}\|_{D(\mathcal{H}^{\beta})} +\sup_{\tau\in
[0,t]}\|y_{2}-z_{2}\|_{D(\mathcal{H}^{\beta})} + \sup_{\tau\in
[0,t]}\|y_{3}-z_{3}\|_{D(\mathcal{H}^{\beta})}\Big)
\end{eqnarray*}
\end{small}
and
\begin{small}
\begin{eqnarray*}
\sup_{t\in [0,T]}\|Ry_{1}(t)-Rz_{1}(t) \|_{D(\mathcal{H}^{\beta})}
&\leq&  \frac{1}{4} \Big( \sup_{\tau\in [0,t]}
\|y_{1}-z_{1}\|_{D(\mathcal{H}^{\beta})} +\sup_{\tau\in
[0,t]}\|y_{2}-z_{2}\|_{D(\mathcal{H}^{\beta})} + \sup_{\tau\in
[0,t]}\|y_{3}-z_{3}\|_{D(\mathcal{H}^{\beta})}\Big).
\end{eqnarray*}
\end{small}
Thus
$$\sup\limits_{t\in [0,T]}\|(P,Q,R)(y_{1}(t),y_{2}(t),y_{3}(t)) -
(P,Q,R)(z_{1}(t),z_{2}(t),z_{3}(t))\|_{D(\mathcal{H}^{\beta})^{3}}$$
\begin{eqnarray*}
 &\leq& \sup_{\tau\in [0,t]} \Big( \|Py_{1}-Pz_{1}\|_{D(\mathcal{H}^{\beta})}
 +\|Qy_{2}-Qz_{2}\|_{D(\mathcal{H}^{\beta})} + \|Ry_{3}-Rz_{3}\|_{D(\mathcal{H}^{\beta})}\Big),\\
&\leq& \sup_{\tau\in [0,t]}
\|Py_{1}-Pz_{1}\|_{D(\mathcal{H}^{\beta})} +\sup_{\tau\in
[0,t]}\|Qy_{2}
-Qz_{2}\|_{D(\mathcal{H}^{\beta})} + \sup_{\tau\in [0,t]}\|Ry_{3}-Rz_{3}\|_{D(\mathcal{H}^{\beta})},\\
&\leq& \frac{3}{4} \Big(\sup_{\tau\in [0,t]}
\|y_{1}-z_{1}\|_{D(\mathcal{H}^{\beta})} +\sup_{\tau\in
[0,t]}\|y_{2}-z_{2}\|_{D(\mathcal{H}^{\beta})} + \sup_{\tau\in
[0,t]}
\|y_{3}-z_{3}\|_{D(\mathcal{H}^{\beta})}\Big),\\
&\leq& \frac{3}{4} \sup_{t\in [0,T]}\|(y_{1}(t),y_{2}(t),y_{3}(t)) -
(z_{1}(t),z_{2}(t),z_{3}(t))\|_{D(\mathcal{H}^{\beta})^{3}}.
\end{eqnarray*}
Hence $(P,Q,R)$ is a strict contraction on
$B_{r_{0}}(H_{0},I_{0},V_{0})$ and this proves Part b). According to
Banach's Fixed Point Theorem, (P,Q,R) has a unique fixed point in
$B_{r_{0}}(H_{0},I_{0},V_{0})$. This is the solution of
(\ref{equation1.1}) on [0,T] with initial value
$(H(0),I(0),V(0))=(H_{0},I_{0},V_{0})$ in
$(D(\mathcal{H}^{\beta}))^{3}$. This completes the proof of
Proposition~\ref{propo3.9}.

\section{Proof of Theorem~\ref{theo2.4.2}}
  We first prove the existence and positivity of the solution.
  From Theorem~\ref{theo2.4.1}, there exist a
sequence $(w^{m})$ and a function $w$ such that
  $$ w^{m}\rightarrow w \;\;\mbox{in}\;\;C([0,T], \mathcal{H}),$$
   with $w\geq 0$ and $w(0)=w_{0}$.\\
   We check that
   $$ q(w^{m-1})\rightarrow q(w)\;\;\mbox{and}\;\;
   q(w^{m-1})w^{m}
  \rightarrow q(w)w\;\;\mbox{in}\;\;C([0,T], \mathcal{H}).$$
 We also have
 $$f(w^{m-1})\rightarrow f(w)\;\; \mbox{in} \;\;C([0,T],
  \mathcal{H})\bigcap L^{2}((0, T), E').$$
  But, $w^{m}$ is solution of
\begin{equation}
  \Big< \frac{\partial w^{m}}{\partial t}, v \Big> +
   \Big<Aw^{m} , v \Big> + \left( q(w^{m-1})w^{m}, v\right)
   =\Big<f(w^{m-1}) , v \Big>,\;\; \forall v\in E.
\end{equation}
  We take $\phi \in \mathcal{D}((0,T))$, so that
  $\phi v \in L^{2}((0,T), E)$,
  \begin{equation}\label{equation2.28}
 \int_{0}^{T} \Big< \frac{\partial w^{m}}{\partial t},
 \phi v \Big> dt + \int_{0}^{T}\Big<Aw^{m} , \phi v \Big>dt
  + \int_{0}^{T}\left( q(w^{m-1})w^{m}, \phi v\right)dt
  =\int_{0}^{T}\Big<f(w^{m-1}) , \phi v \Big>dt.
\end{equation}
The second term in the left side and the right side of the equality
(\ref{equation2.28}) converges due to the weak convergence in
$L^{2}((0, T), E')$. The third term in the left-hand side of
(\ref{equation2.28}) also converges, due to the convergence in
$C([0,T], \mathcal{H})$. We deduce that $\frac{\partial
w^{m}}{\partial t}$ converges weakly in $L^{2}((0, T), E')$.
\\\indent
But we have
$$w^{m}\rightarrow w \;\;\mbox{in}\;\; C([0,T], \mathcal{H}).$$
Then
\begin{equation*}
  \frac{\partial w^{m}}{\partial t} \rightarrow
 \frac{\partial w}{\partial t}\;\; \mbox{in}\;\; \mathcal{D'}((0, T),
 \mathcal{H})
\end{equation*}
Therefore, we obtain
$$\frac{\partial w^{m}}{\partial t} \rightarrow
\frac{\partial w}{\partial t}\;\;\mbox{weakly in}\;\; L^{2}((0, T),
E'),$$ and
  \begin{equation}\label{equation2.29}
 \int_{0}^{T} \Big< \frac{\partial w}{\partial t},
 \phi v \Big> dt + \int_{0}^{T}\Big<Aw , \phi v
 \Big>dt + \int_{0}^{T}\Big( q(w)w, \phi v\Big)_
 {\mathcal{H}}dt=\int_{0}^{T}\Big<f(w) , \phi v \Big>dt.
\end{equation}
This being true for all  $\phi$, one has
  $$ \Big<
\frac{\partial w}{\partial t}, v \Big>  + \Big<Aw , v \Big> + \Big(
q(w)w, v\Big) _{\mathcal{H}}=\Big<f(w) , v \Big>, \;\; \forall v \in
E .$$ That is to say,
  \begin{eqnarray}\label{equation2.30}
  \frac{d}{d t}(w,v)_{\mathcal{H}} + a(w , v)+
  \Big( q(w)w, v\Big)_{\mathcal{H}}=\Big<f(w) , v \Big>,\;
   \forall v \in E,\\
  \frac{\partial w}{\partial t}=f(w)-Aw-q(w)w
  \;\;\mbox{in}\;\; L^{2}((0, T), E').
\end{eqnarray}
According to (\ref{equation2.023}) and (\ref{equation2.024}) we have
\begin{eqnarray}\label{equation2.31}
 w^{m}(t)&=& G(t)w_{0} + \int_{0}^{t} G(t-s)(-q(w^{m-1})w^{m}
  + f(w^{m-1}))(s)ds,
\end{eqnarray}
and in addition, as $q(w^{m-1})w^{m}$ and $f(w^{m-1})$ converge
 in
$C^{0}([0,T], \mathcal{H})$ and the operator $\mathcal{G}$, defined
by the relation (\ref{equation2.025}), is compact, using the limit
in (\ref{equation2.31}) one has,
\begin{eqnarray}\label{equation2.32}
 w(t)&=& G(t)w_{0} + \int_{0}^{t} G(t-s)(-q(w)w + f(w))(s)ds.
\end{eqnarray}
It remains to prove uniqueness.
\\\indent
Let $v$ be another solution of IBVP (\ref{equation2.15}). Then
 $$
v\in W(0, T, E, E') \Rightarrow v\in C([0,T], \mathcal{H})\;\;
 \mbox{and}\;\; v\geq0.$$
Consequently we obtain
 $$q(v)v+f(v) \in L^{2}((0, T), E').$$
Thus, by Proposition~2.11 of \cite{WS}, one has
\begin{eqnarray*}
 v(t)&=& G(t)w_{0} + \int_{0}^{t} G(t-s)(-q(v)v + f(v))(s)ds.
\end{eqnarray*}
Subtracting, we have
\begin{eqnarray}\label{equation2.33}
 w(t)-v(t)= \int_{0}^{t} G(t-s)\Big(-(q(w)w-q(v)v) + (f(w)-f(v))
 \Big)(s)ds,
\end{eqnarray}
with
\begin{eqnarray}
q(w)w-q(v)v &=& q(w)w-q(w)v + q(w)v -q(v)v, \nonumber\\
             &=& q(w)(w-v) + (q(w)-q(v))v. \nonumber
\end{eqnarray}
Since $w$ is positive, one has
\begin{eqnarray*}
\Bigg\|
\frac{w_{j}}{\alpha_{0}+\alpha_{1}w_{k}+\alpha_{2}w_{j}+\alpha_{3}w_{k}w_{j}}
\Bigg\| \leq \frac{1}{K} \|w_{j}\|_{\infty}
\end{eqnarray*}
where
\begin{eqnarray*}
\|w_{j}\|_{\infty}=\|w_{j}\|_{L^{\infty}((0, T), \mathcal{H})}.
\end{eqnarray*}
If we define
$$\|w\|_{\infty}= \sum_{j=1}^{3}\|w_{j}\|_{\infty},$$
there is $M_{1}>0$ such that
 \begin{eqnarray*}
 \|q(w)\|_{\infty}\leq M_{1}\|w_{j}\|_{\infty}.
 \end{eqnarray*}
So, for $r=1,2,3$, the numerator of  $q_{r}(w)-q_{r}(v)$ is the sum
of terms of the form $(w_{k}-v_{k})v_{j}$ or $(w_{j}-v_{j})w_{k}$,
and we can find $M_{2}>0$ such that
\begin{eqnarray*}
\big| q_{r}(w)-q_{r}(v) \big |_{\mathcal{H}}(s) \leq
M_{2}\Big(\sum_{j=1}^{3}|w_{j}(s)-v_{j}(s)|_{\mathcal{H}}\Big).
\end{eqnarray*}
Also we can find $M_{3}>0$ such that
\begin{equation*}
\big| f_{r}(w)-f_{r}(v) \big |_{\mathcal{H}}(s) \leq
M_{3}\Big(\sum_{j=1}^{3}|w_{j}(s)-v_{j}(s)|_{\mathcal{H}}\Big).
\end{equation*}
Summing up  $|w_{j}(s)-v_{j}(s)|_{\mathcal{H}}$ and noting that
$\|G_{j}(t-s)\|\leq N_{j}e^{\theta_{j}T}$ with
$N_{j},\;\theta_{j}>0$, we can find $M>0$ such that
\begin{eqnarray*}
\sum_{j=1}^{3}|w_{j}(s)-v_{j}(s)|_{\mathcal{H}}\leq
M\|w-v\|_{\infty}.
\end{eqnarray*}
Replacing in (\ref{equation2.33}), we obtain
\begin{eqnarray*}
\sum_{j=1}^{3}|w_{j}(s)-v_{j}(s)|_{\mathcal{H}}\leq
M^{2}\|w-v\|_{\infty}\int_{0}^{t}sds=M^{2}\frac{t^{2}}{2}
\|w-v\|_{\infty}.
\end{eqnarray*}
By induction, we have
\begin{equation*}
\sum_{j=1}^{3}|w_{j}(s)-v_{j}(s)|_{\mathcal{H}}\leq
 \frac{M^{n}}{n!}T^{n}\|w-v\|_{\infty},
\end{equation*}
with
\begin{equation*}
\lim_{n\rightarrow +\infty}\frac{M^{n}}{n!}T^{n}\|w-v\|_{\infty}=0.
\end{equation*}
Therefore $w=v$. This ends the proof of Theorem~\ref{theo2.4.2}.

\end{document}